\documentclass[11pt]{article}

\usepackage[a4paper,margin=0.92in,headheight=15pt]{geometry}
\usepackage{amsmath,amssymb,amsthm,mathtools,mathrsfs}
\usepackage{booktabs,array,tabularx,longtable,multirow}
\usepackage{enumitem}
\usepackage{microtype}
\usepackage{graphicx}
\usepackage{caption,subcaption,float}
\usepackage{pdflscape}
\usepackage[numbers,sort&compress]{natbib}
\usepackage{xcolor}
\usepackage{fancyhdr}
\usepackage{url}
\usepackage[hidelinks]{hyperref}
\usepackage[nameinlink,noabbrev]{cleveref}

\graphicspath{{figures/}}

\definecolor{deepblue}{RGB}{21,56,100}
\definecolor{softgray}{RGB}{242,244,247}
\definecolor{darkgray}{RGB}{70,70,70}

\hypersetup{
  pdftitle={A Topological Perspective on the Birch and Swinnerton-Dyer Conjecture},
  pdfauthor={Maisara Shoeib},
  pdfsubject={Elliptic curves, Mordell-Weil lattices, flat tori, BSD conjecture},
  pdfkeywords={Birch and Swinnerton-Dyer conjecture, elliptic curves, Mordell-Weil lattice, canonical height, regulator, flat torus, theta series}
}

\newtheorem{theorem}{Theorem}[section]
\newtheorem{proposition}[theorem]{Proposition}
\newtheorem{corollary}[theorem]{Corollary}
\newtheorem{lemma}[theorem]{Lemma}

\theoremstyle{definition}
\newtheorem{definition}[theorem]{Definition}
\newtheorem{problem}[theorem]{Bridge Problem}
\newtheorem{protocol}[theorem]{Computational Protocol}
\theoremstyle{remark}
\newtheorem{remark}[theorem]{Remark}
\newtheorem{warning}[theorem]{Warning}

\newcommand{\Q}{\mathbb{Q}}
\newcommand{\R}{\mathbb{R}}
\newcommand{\Z}{\mathbb{Z}}

\newcommand{\F}{\mathbb{F}}
\newcommand{\T}{\mathcal{T}}

\newcommand{\Vol}{\operatorname{Vol}}
\newcommand{\Reg}{\operatorname{Reg}}
\newcommand{\rank}{\operatorname{rank}}

\newcommand{\sys}{\operatorname{sys}}
\newcommand{\Spec}{\operatorname{Spec}}
\newcommand{\Tr}{\operatorname{Tr}}
\newcommand{\ord}{\operatorname{ord}}
\newcommand{\Hom}{\operatorname{Hom}}

\newcommand{\Sha}{\mathop{\mathrm{Sha}}\nolimits}

\newcommand{\CP}{\mathrm{CP}}
\newcommand{\Sil}{\mathrm{Sil}}
\newcommand{\e}{\mathrm{e}}
\newcommand{\dd}{\,\mathrm{d}}
\newcommand{\abs}[1]{\left\lvert#1\right\rvert}
\newcommand{\norm}[1]{\left\lVert#1\right\rVert}
\newcommand{\inner}[2]{\left\langle#1,#2\right\rangle_E}
\newcommand{\qheight}{\widehat h_E^{\CP}}

\setlist[itemize]{leftmargin=1.4em,itemsep=0.25em,topsep=0.35em}
\setlist[enumerate]{leftmargin=1.8em,itemsep=0.3em,topsep=0.4em}

\title{\textbf{A Topological Perspective on the Birch and Swinnerton--Dyer Conjecture}}
\author{Maisara Shoeib\\[2pt]
\small Higher Colleges of Technology (HCT), United Arab Emirates\\
\small \texttt{mshoeib@hct.ac.ae}}
\date{}

\begin{document}
\maketitle

\begin{abstract}
For an elliptic curve $E/\Q$, let $M_E=E(\Q)/E(\Q)_{\mathrm{tors}}$, let $V_E=M_E\otimes_{\Z}\R$, and equip $V_E$ with the positive-definite N\'eron--Tate pairing in the normalization used in the computational BSD formula.  The quotient $\T_E=V_E/M_E$ is a compact flat torus, called the Mordell--Weil height torus.  We prove that $H_1(\T_E,\Z)\cong M_E$, $b_1(\T_E)=\rank E(\Q)$, the marked closed geodesic attached to $P\in M_E$ has squared length $\widehat h_E(P)$, and $\Vol(\T_E)^2=\Reg(E/\Q)$.  We also derive marked-length reconstruction, systolic and successive-minimum bounds, finite-index covering laws, isogeny scaling, toroidal-helical representatives, four-dimensional character windows, theta identities, and exact heat-trace formulae.

Independently, we construct an explicit-formula bilinear form from the completed $L$-function.  For band-limited Paley--Wiener seeds $f,g$, the test function $h_{f,g}=fg$ has compactly supported Fourier transform, so the prime-power contribution is exactly finite.  Comparing this analytic form with a positive Hodge-trace form on $\T_E$ defines the height-torus residual.  Its central-mode decomposition contains the term
\[
  \bigl(\ord_{s=1}L(E,s)-\rank E(\Q)\bigr)f(0)g(0),
\]
while the remaining residual measures the noncentral spectral mismatch.  We prove congruence covariance, unitary-phase invariance, and a counting-law obstruction: the raw compact height-torus spectrum grows as $\Theta(T^r)$, whereas the zeros of the degree-two elliptic-curve $L$-function grow as $\Theta(T\log T)$, so the raw torus spectrum cannot be the complete zero spectrum under a fixed affine rescaling.

Thus the rank part of BSD is exactly reformulated as $\ord_{s=1}L(E,s)=b_1(\T_E)$, and the regulator in the leading-coefficient formula becomes $\Vol(\T_E)^2$.  These are reformulations, not a proof of BSD.  Computations for ranks $0$--$3$ verify the internal height-torus dictionary, database-basis status, bounded saturation checks, normalization, and BSD consistency quotients; no unvalidated residual-sign claim is made.
\end{abstract}

\noindent\textbf{Keywords.} Elliptic curves; Birch and Swinnerton--Dyer conjecture; Mordell--Weil lattice; N\'eron--Tate height; flat torus; regulator; closed geodesics; theta series; heat trace; explicit formula; spectral residual; spectral geometry.

\medskip
\noindent\textbf{2020 Mathematics Subject Classification.} 11G40, 11G05, 14G05, 14H52, 53C22, 58J50, 11H06.

\section{Introduction}

The Birch and Swinnerton--Dyer conjecture (BSD) connects two invariants that arise in apparently different mathematical languages.  The algebraic rank of $E(\Q)$ is defined by rational points and the group law, whereas the analytic rank is the order of vanishing of the Hasse--Weil $L$-function at its central point.  The original numerical investigations of Birch and Swinnerton-Dyer \cite{BirchSwinnertonDyer1963,BirchSwinnertonDyer1965} suggested the equality
\begin{equation}\label{eq:BSD-rank-intro}
  \rank E(\Q)=\ord_{s=1}L(E,s).
\end{equation}
The full conjecture predicts the first nonzero Taylor coefficient in terms of the real period, the N\'eron--Tate regulator, the Tate--Shafarevich group, the Tamagawa factors, and the torsion subgroup \cite{Cassels1965,Tate1966,WilesClay2006}.

The geometric intuition behind the present work is that a free abelian group of rank $r$ should admit $r$ independent periodic directions.  The important issue is to identify the correct space on which those directions become genuine loops.  The real or complex elliptic curve itself is not the correct space: its topological type is fixed independently of the arithmetic rank.  Nor does the countable orbit $\{nP:n\in\Z\}$ become a continuous loop merely by applying bounded coordinate functions.  The correct object is instead obtained from the free Mordell--Weil group as a Euclidean lattice under the canonical height.

The theory of Mordell--Weil lattices is well established \cite{Neron1965,Silverman2009,Shioda1990,Shioda1991,SchuettShioda2019}.  High-rank constructions and extremal examples further illustrate how rich the associated height lattices can be \cite{Elkies2000,Elkies2007}.  The point of this paper is not to rename a known lattice and then claim a solution of BSD.  It is to assemble the arithmetic lattice, its quotient flat torus, its marked geodesics, its volume, and its heat spectrum into one exact dictionary tailored to BSD; to correct the specific topological obstructions that arise in naive orbit embeddings; to state several useful consequences explicitly; and to isolate the genuinely open analytic step without importing it into a definition.

\begin{figure}[H]
  \centering
  \includegraphics[width=0.98\textwidth]{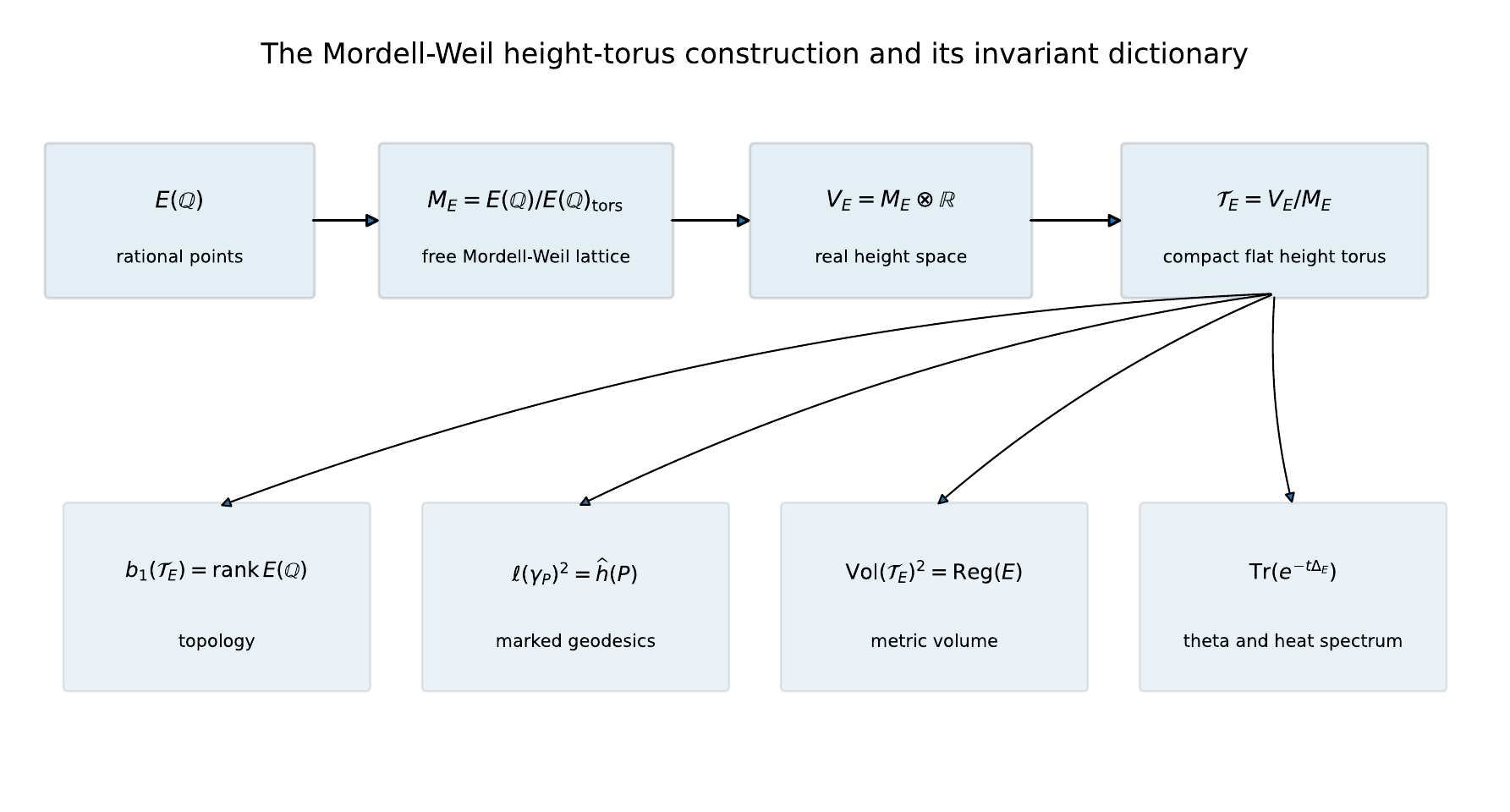}
  \caption{The construction $E(\Q)\to M_E\to V_E\to\T_E$ and the resulting invariant dictionary.  The topology records the rank, the marked geodesic lengths record the height form, the volume records the regulator, and the heat trace records the dual lattice spectrum.  The independent analytic comparison is introduced in \cref{sec:analytic-side}.}
  \label{fig:pipeline}
\end{figure}

\subsection{Principal results}

Let
\[
  M_E:=E(\Q)/E(\Q)_{\mathrm{tors}},
  \qquad
  V_E:=M_E\otimes_{\Z}\R,
\]
and equip $V_E$ with the BSD-normalized N\'eron--Tate pairing.  The central object is
\[
  \T_E:=V_E/M_E.
\]
The paper establishes the following statements.

\begin{enumerate}[label=\textup{(\arabic*)}]
  \item \textbf{Rank and homology.}
  \[
    \pi_1(\T_E)\cong H_1(\T_E,\Z)\cong M_E,
    \qquad
    b_1(\T_E)=\rank E(\Q).
  \]
  More generally, $H_k(\T_E,\Z)\cong\bigwedge^k M_E$.

  \item \textbf{Height and marked geodesics.}  The class $P\in M_E$ determines a genuine closed geodesic $\gamma_P(t)=tP+M_E$, and
  \[
    \operatorname{Length}(\gamma_P)^2=\qheight(P).
  \]
  A finite family of marked lengths reconstructs the entire height Gram matrix.

  \item \textbf{Regulator and volume.}
  \[
    \Vol(\T_E)^2=\Reg(E/\Q).
  \]
  The systole and successive minima of $\T_E$ are the square roots of the corresponding minimal canonical heights.

  \item \textbf{Finite-index control.}  An unsaturated rank-$r$ subgroup produces a finite covering torus whose volume is inflated by the index and whose Gram determinant is inflated by the square of the index.  This gives a geometric interpretation of saturation errors in computation.

  \item \textbf{Helical models and four-dimensional windows.}  Every Mordell--Weil class has a toroidal-helical representative with an integral winding vector.  A single four-dimensional character window can retain at most two independent homology directions; $\lceil r/2\rceil$ such windows are necessary and sufficient within this construction.

  \item \textbf{Spectral formulation.}  The scalar Laplace spectrum of $\T_E$ is the squared-norm spectrum of the dual Mordell--Weil lattice.  Poisson summation gives an exact heat-trace formula from which the rank, volume, and systole can be recovered asymptotically.

  \item \textbf{Topological BSD.}  The rank conjecture is exactly equivalent to
  \[
    \ord_{s=1}L(E,s)=b_1(\T_E)
    =\dim\ker\Delta_E^{(1)},
  \]
  and the regulator in the full BSD formula is $\Vol(\T_E)^2$.

  \item \textbf{Independent analytic form and exact local support.}  A centered explicit-formula functional produces $G^{\mathrm{an}}_{E,\Delta}$ without using a Mordell--Weil basis.  For band-limited Paley--Wiener seeds $f,g$, the product $h_{f,g}=fg$ has $\operatorname{supp}\widehat h_{f,g}\subset[-\Delta,\Delta]$, so only prime powers $n\le \exp(2\pi\Delta)$ occur.

  \item \textbf{Central-mode residual and obstruction.}  The analytic--torus residual is basis-covariant and decomposes as
  \[
    \mathcal R_{E,\Delta,\tau}(f,g)
    =\bigl(r_{\mathrm{an}}(E)-r\bigr)f(0)g(0)
     +\mathcal R^{\circ}_{E,\Delta,\tau}(f,g).
  \]
  The raw torus spectrum cannot equal the complete zero spectrum because their counting laws are respectively $\Theta(T^r)$ and $\Theta(T\log T)$.

  \item \textbf{Certificate layer.}  A finite residual identity, one-sided comparison, controlled norm bound, or obstruction is accepted only after independent construction, interval or exact error control, fixed acceptance margins, and negative controls.
\end{enumerate}

\subsection{Status and non-claims}

The distinction between theorem, reformulation, and open problem is essential.

\begin{itemize}
  \item The lattice, torus, homology, length, volume, theta, heat, and isogeny statements proved below are rigorous consequences of the Mordell--Weil theorem, canonical heights, and standard flat-torus geometry.
  \item The equality $\ord_{s=1}L(E,s)=b_1(\T_E)$ is exactly BSD rewritten through a proved identity $b_1(\T_E)=\rank E(\Q)$.  Rewriting does not prove the analytic equality.
  \item The helical drawings are coordinate representations of homology classes on $\T_E$.  They are not algebraic parametrizations of $E$, and their Euclidean appearance depends on a chosen basis and projection.
  \item No operator is declared to have the zeros of $L(E,s)$ as its spectrum.  Such a declaration would be circular unless the operator, its domain, self-adjointness, spectral determinant, and relation to the Euler product were independently established.
  \item The explicit-formula matrix and the height-torus trace matrix are defined by separate pipelines.  Their residual is not assumed to vanish or have a preferred sign.
  \item The spectral-counting obstruction proved below rules out the raw compact torus as the complete zero operator, while leaving the central-multiplicity formulation of BSD open.
\end{itemize}

Classical progress includes the complex-multiplication cases of Coates--Wiles, the Gross--Zagier and Gross--Kohnen--Zagier height formulae, Kolyvagin's Euler systems, and Iwasawa-theoretic main-conjecture results \cite{CoatesWiles1977,GrossZagier1986,GrossKohnenZagier1987,Kolyvagin1988,Rubin1991}.  Recent work continues to extend the range in which parts of BSD are known, including Iwasawa-theoretic results and infinite twist families \cite{SkinnerUrban2014,JetchevSkinnerWan2017,BurungaleSkinnerTianWan2024,BanwaitHuang2026}.  Statistical, average-rank, and family-level results provide further evidence and constraints on possible rank distributions \cite{HeathBrown2004,BhargavaShankarBinary2015,BhargavaShankarTernary2015,BhargavaSkinner2014,BhargavaSkinnerZhang2014,ParkPoonenVoightWood2019,Smith2025}.  None of these developments removes the need for a non-circular bridge in general rank.

A distinct but directly relevant research tradition seeks cohomological and dynamical realizations of zeta and $L$-functions, together with realizations in noncommutative geometry.  Deninger formulated regularized-determinant and dynamical analogies for local factors and arithmetic zeta functions \cite{Deninger1991,Deninger1992,Deninger2002}; Connes interpreted the explicit formula as a trace formula on the noncommutative ad\`ele-class space \cite{Connes1999}; and Hesselholt realized a Deninger-type regularized-determinant description of Hasse--Weil zeta functions for smooth proper schemes over finite fields using topological Hochschild homology \cite{Hesselholt2016}.  These programs are the closest established antecedents to the bridge problem formulated below.  The Mordell--Weil height torus constructed here is not claimed to realize their arithmetic cohomology or trace operators: it supplies the rigorously defined target-side geometry that any elliptic-curve bridge would still have to connect to the independently defined $L$-function.

\section{Arithmetic background and normalization}

\subsection{Elliptic curves and the Mordell--Weil group}

Let $E/\Q$ be an elliptic curve written in a global Weierstrass form
\begin{equation}\label{eq:general-Weierstrass}
  y^2+a_1xy+a_3y=x^3+a_2x^2+a_4x+a_6,
\end{equation}
with nonzero discriminant.  The Mordell--Weil theorem \cite{Mordell1922,Weil1929}, in the standard elliptic-curve formulation of \cite{Cassels1991,Silverman2009}, states that
\begin{equation}\label{eq:MW-decomp}
  E(\Q)\cong E(\Q)_{\mathrm{tors}}\oplus\Z^r.
\end{equation}
We set
\[
  M_E:=E(\Q)/E(\Q)_{\mathrm{tors}},
\]
so $M_E$ is a free abelian group of rank $r$.

\subsection{Cremona/PARI canonical-height normalization}

Canonical-height conventions in the literature differ by a factor of two.  This paper uses throughout the normalization implemented by PARI/GP's \texttt{ellheight} and \texttt{ellheightmatrix}, which follows Cremona's \emph{Algorithms for Modular Elliptic Curves} \cite{Cremona1997,PARI2026}.  We call it the Cremona/PARI normalization and denote it by $\widehat h_E^{\CP}$.  This is the normalization whose full-basis Gram determinant is the regulator paired with \texttt{ellbsd} in the numerical BSD consistency quotient.

\begin{definition}[Cremona/PARI height form]\label{def:height-normalization}
For $P,Q\in E(\Q)$, let
\[
  \qheight(P):=\widehat h_E^{\CP}(P),
\]
and define the associated bilinear form by either equivalent polarization identity
\begin{equation}\label{eq:height-pairing}
  \inner{P}{Q}
  :=\frac14\bigl(\qheight(P+Q)-\qheight(P-Q)\bigr)
  =\frac12\bigl(\qheight(P+Q)-\qheight(P)-\qheight(Q)\bigr).
\end{equation}
Then $\inner{P}{P}=\qheight(P)$.
\end{definition}

The form descends to $M_E$ and is positive definite there \cite{Neron1965,Silverman2009,HindrySilverman2000}.  The PARI/GP reference manual states explicitly that its height matrix should be divided by $2$ to agree with the normalization used, for example, in Silverman's books; in symbols,
\begin{equation}\label{eq:cp-silverman-conversion}
  \widehat h_E^{\CP}=2\widehat h_E^{\Sil}.
\end{equation}
Consequently, conversion to the Silverman convention divides Gram entries and squared geodesic lengths by $2$, the rank-$r$ regulator by $2^r$, and the volume by $2^{r/2}$.  No such division is made in the tables or in the \texttt{ellbsd} quotients below: those calculations use the direct Cremona/PARI output consistently.  Dividing only the regulator while retaining the unconverted \texttt{ellbsd} factor would mix conventions and produce a spurious factor of $2^r$.

For a basis $\mathcal B=(P_1,\ldots,P_r)$ of $M_E$, define the Gram matrix
\begin{equation}\label{eq:Gram}
  G_E(\mathcal B):=\bigl(\inner{P_i}{P_j}\bigr)_{1\le i,j\le r}.
\end{equation}
The Cremona/PARI regulator used in this manuscript is
\begin{equation}\label{eq:regulator}
  \Reg(E/\Q):=\det G_E(\mathcal B),
\end{equation}
with the convention $\Reg(E/\Q)=1$ when $r=0$.  Since a basis change is given by $A\in\mathrm{GL}_r(\Z)$ and sends $G$ to $A^{\mathsf T}GA$, the determinant is basis independent.  All metric theorems remain valid in any other height convention after the corresponding uniform rescaling.

\subsection{The Hasse--Weil \texorpdfstring{$L$}{L}-function and BSD}

For each prime $p$ of good reduction, let
\[
  a_p=p+1-\#E(\F_p).
\]
Hasse's bound gives $\abs{a_p}\le2\sqrt p$.  The $L$-function is represented in its half-plane of absolute convergence by the Euler product
\begin{equation}\label{eq:L-function}
  L(E,s)=
  \prod_{p\nmid N_E}
  \bigl(1-a_pp^{-s}+p^{1-2s}\bigr)^{-1}
  \prod_{p\mid N_E}
  \bigl(1-a_pp^{-s}\bigr)^{-1}.
\end{equation}
Modularity \cite{Wiles1995,TaylorWiles1995,BreuilConradDiamondTaylor2001} supplies analytic continuation and a functional equation centered at $s=1$.

Set
\[
  r_{\mathrm{an}}(E):=\ord_{s=1}L(E,s),
  \qquad
  L^*(E,1):=\lim_{s\to1}\frac{L(E,s)}{(s-1)^{r_{\mathrm{an}}(E)}}.
\]
The rank part of BSD is $r_{\mathrm{an}}(E)=r$.  The leading-coefficient formula predicts, with compatible standard choices of period and height normalization,
\begin{equation}\label{eq:BSD-full}
  L^*(E,1)
  =
  \frac{
    \Omega_E\,\Reg(E/\Q)\,\#\Sha(E/\Q)\,\prod_p c_p
  }{
    \#E(\Q)_{\mathrm{tors}}^2
  }.
\end{equation}
The rank equality and finiteness of $\Sha$ are known when the analytic rank is at most one through Gross--Zagier and Kolyvagin, together with modularity \cite{GrossZagier1986,Kolyvagin1988}; the complex-multiplication and Iwasawa-theoretic foundations include \cite{CoatesWiles1977,Rubin1991}.  The full leading-coefficient formula involves additional local and $p$-adic work and is known in important but not universal generality \cite{Kato2004,SkinnerUrban2014,JetchevSkinnerWan2017,BurungaleSkinnerTianWan2024}.

For analytic rank at least two, a numerical value obtained by solving \eqref{eq:BSD-full} for $\#\Sha$ is only the \emph{BSD-predicted} order of $\Sha$.  Without an independent theorem proving finiteness and controlling the relevant Selmer groups, an apparently integral value---even the value $1$ to high precision---is a consistency check, not a proof that $\Sha$ is finite or trivial.

\section{Why the naive real-locus and discrete-orbit models fail}

A rigorous topological framework must first identify the obstruction in the most immediate geometric constructions.  This section records three such obstructions.  They do not reject the topological objective; they determine the space on which it must be carried out.

\subsection{The real locus has bounded topological complexity}

\begin{theorem}[Real-locus obstruction]\label{thm:real-locus-obstruction}
For an elliptic curve $E/\R$, the compact one-dimensional Lie group $E(\R)$ has either one or two connected components, each homeomorphic to $S^1$.  Consequently,
\[
  b_1(E(\R))\in\{1,2\}.
\]
If $\Phi:E(\R)\hookrightarrow\R^n$ is a topological embedding, then
\[
  b_1(\Phi(E(\R)))=b_1(E(\R))\le2.
\]
Therefore an embedding of the real elliptic curve cannot encode arbitrary Mordell--Weil rank by the number of independent loops in its image.
\end{theorem}

\begin{proof}
The real points of a smooth projective cubic form a compact one-dimensional real Lie group.  Its identity component is a compact connected one-dimensional Lie group and is therefore a circle.  The component group of a real elliptic curve has order one or two.  Thus $E(\R)$ is one circle or a disjoint union of two circles.  An embedding is a homeomorphism onto its image and preserves singular homology.
\end{proof}

The obstruction is decisive: curves of algebraic rank $3$, $8$, or higher exist, while the topology of $E(\R)$ remains one or two circles.  The arithmetic rank must therefore be represented in a space built from the arithmetic group, not inferred from the unmodified real curve.

\subsection{A countable orbit is not a continuous loop}

\begin{proposition}[Discrete-orbit obstruction]\label{prop:discrete-orbit-obstruction}
Let $X$ be a metric space, let $P$ be a point of infinite order in an abelian group, and let $\Phi$ be any map from the group to $X$.  The orbit set
\[
  \mathcal O_P:=\{\Phi(nP):n\in\Z\}
\]
is countable.  It cannot be the image of a nonconstant continuous loop whose image is connected and contained in $\mathcal O_P$.
\end{proposition}

\begin{proof}
The continuous image of $S^1$ is connected.  Every connected metric space containing two distinct points is uncountable.  Since $\mathcal O_P$ is countable, every connected subset of $\mathcal O_P$ that is the image of a loop must be a singleton.
\end{proof}

The closure of an orbit may be a circle or a higher-dimensional subtorus, but the orbit, its closure, and a parametrized continuous loop are different objects.  Homology classes require an actual continuous singular cycle.  The construction in \cref{sec:height-torus} supplies one canonically.

\subsection{A boundedness obstruction for the former prime sum}

Consider the statistic
\begin{equation}\label{eq:F-old}
  F_{\mathrm{old}}(E,X)
  :=\frac1X\sum_{\substack{p\le X\\p\nmid N_E}}
  \frac{a_p\log p}{\sqrt p}.
\end{equation}

\begin{proposition}[Hasse-bound obstruction]\label{prop:F-old-bounded}
For every elliptic curve $E/\Q$,
\[
  F_{\mathrm{old}}(E,X)=O(1).
\]
More explicitly,
\[
  \abs{F_{\mathrm{old}}(E,X)}
  \le \frac{2}{X}\sum_{p\le X}\log p.
\]
Consequently, $F_{\mathrm{old}}(E,X)$ cannot have asymptotic growth $C(\log X)^r$ with $C\ne0$ and $r>0$.
\end{proposition}

\begin{proof}
Hasse's inequality gives $\abs{a_p}/\sqrt p\le2$.  Hence
\[
  \abs{F_{\mathrm{old}}(E,X)}
  \le \frac1X\sum_{p\le X}\frac{\abs{a_p}\log p}{\sqrt p}
  \le \frac{2}{X}\sum_{p\le X}\log p.
\]
The Chebyshev function $\vartheta(X)=\sum_{p\le X}\log p$ is $O(X)$, proving boundedness.
\end{proof}

This does not imply that all finite prime sums are useless.  Properly normalized Mestre--Nagao-type statistics can correlate with analytic rank, and recent work has found conductor-dependent oscillations and improvements from multivalue or learned weights \cite{Nagao1997,RosenSilverman1998,KazalickiVlah2022,BujanovicKazalickiNovak2024,BujanovicKazalickiVlah2025,BieriEtAl2026}.  It does imply that the particular power-growth law in \eqref{eq:F-old} must not be used.

\section{The Mordell--Weil height torus}\label{sec:height-torus}

\begin{definition}[Height space and height torus]\label{def:height-torus}
Let
\[
  V_E:=M_E\otimes_{\Z}\R,
\]
and extend the pairing \eqref{eq:height-pairing} by $\R$-bilinearity.  The \emph{Mordell--Weil height torus} is
\begin{equation}\label{eq:height-torus}
  \T_E:=V_E/M_E
\end{equation}
with the quotient flat Riemannian metric.  When $r=0$, $V_E=0$ and $\T_E$ is a point of volume one.
\end{definition}

\begin{proposition}[Basis independence]\label{prop:basis-independence}
The metric torus $\T_E$ is intrinsic to the paired lattice $(M_E,\inner{\cdot}{\cdot})$ and is independent, up to canonical isometry, of the chosen Mordell--Weil basis.
\end{proposition}

\begin{proof}
A basis identifies $V_E$ with $\R^r$ and $M_E$ with $\Z^r$, with metric matrix $G_E$.  Replacing the basis by $A\in\mathrm{GL}_r(\Z)$ replaces $G_E$ by $A^{\mathsf T}G_EA$ and the lattice coordinates by $A^{-1}\Z^r$.  The induced linear map descends to an isometry of the quotient tori.
\end{proof}

\begin{theorem}[Rank--homology theorem]\label{thm:rank-homology}
There are natural isomorphisms
\begin{equation}\label{eq:fundamental-homology}
  \pi_1(\T_E)\cong M_E,
  \qquad
  H_1(\T_E,\Z)\cong M_E.
\end{equation}
Consequently,
\begin{equation}\label{eq:b1-rank}
  b_1(\T_E)=\rank E(\Q).
\end{equation}
\end{theorem}

\begin{proof}
Choose a basis $P_1,\ldots,P_r$ of $M_E$.  Then $V_E\cong\R^r$, $M_E\cong\Z^r$, and $\T_E\cong\R^r/\Z^r$.  The universal cover is the contractible space $V_E$, and the deck transformations are translations by elements of $M_E$.  Thus $\pi_1(\T_E)\cong M_E$.  The first homology is the abelianization of the fundamental group; $M_E$ is already abelian.
\end{proof}

\begin{theorem}[Full homology and cohomology]\label{thm:full-homology}
For $0\le k\le r$,
\begin{equation}\label{eq:full-homology}
  H_k(\T_E,\Z)\cong\bigwedge^k M_E,
  \qquad
  H^k(\T_E,\Z)\cong\bigwedge^k\Hom(M_E,\Z).
\end{equation}
Hence
\begin{equation}\label{eq:betti}
  b_k(\T_E)=\binom{r}{k},
  \qquad
  P_{\T_E}(t):=\sum_{k=0}^r b_kt^k=(1+t)^r.
\end{equation}
In particular, $\chi(\T_E)=0$ for $r>0$, while $\chi(\T_E)=1$ for $r=0$.
\end{theorem}

\begin{proof}
After a basis choice, $\T_E\cong(S^1)^r$.  The K\"unneth theorem identifies its integral homology algebra with the exterior algebra on $H_1$ \cite{Hatcher2002}.  The cohomology statement follows by the universal coefficient theorem.
\end{proof}

\begin{figure}[H]
  \centering
  \includegraphics[width=0.93\textwidth]{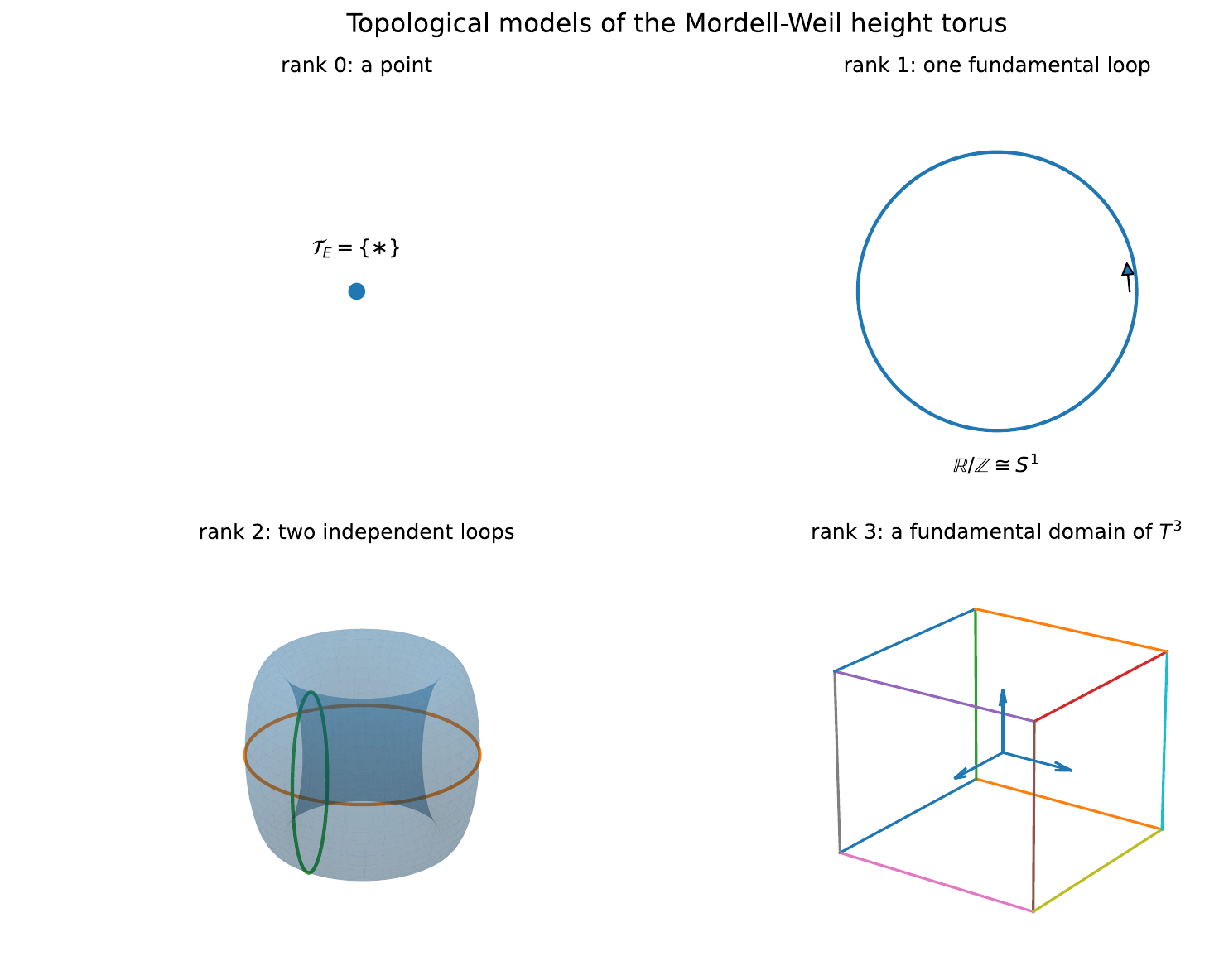}
  \caption{Topological models of $\T_E$ for ranks $0$--$3$.  The rank-three panel shows a fundamental domain of $\R^3/\Z^3$ rather than a globally embedded three-torus in $\R^3$.}
  \label{fig:rank-atlas}
\end{figure}

\begin{remark}[What the torus remembers and forgets]\label{rem:torsion}
The torus records the free Mordell--Weil group and its height metric.  It does not record the torsion subgroup, local Tamagawa factors, the real period, or $\Sha$ by itself.  Those quantities remain separate terms in the full BSD formula.  This limitation is structural, not a defect to be hidden.
\end{remark}

\section{Rational points as marked closed geodesics}

For $P\in M_E$, define
\begin{equation}\label{eq:gamma-P}
  \gamma_P:[0,1]\longrightarrow\T_E,
  \qquad
  \gamma_P(t):=tP+M_E.
\end{equation}
The endpoints agree because $P\in M_E$.

\begin{theorem}[Geodesic realization and height]\label{thm:height-length}
The map
\begin{equation}\label{eq:marked-isomorphism}
  M_E\longrightarrow H_1(\T_E,\Z),
  \qquad
  P\longmapsto[\gamma_P],
\end{equation}
is the natural isomorphism in \cref{thm:rank-homology}.  The loop $\gamma_P$ is a closed geodesic, it minimizes length in its free homotopy class, and
\begin{equation}\label{eq:length-height}
  \operatorname{Length}(\gamma_P)^2=\qheight(P).
\end{equation}
Thus the stable norm on $H_1(\T_E,\R)$ is exactly the N\'eron--Tate norm.
\end{theorem}

\begin{proof}
The lift of $\gamma_P$ to $V_E$ is the straight segment $t\mapsto tP$, whose endpoint displacement is the deck transformation $P$.  Therefore its homotopy and homology class is $P$.  A straight segment is the shortest path between its endpoints in Euclidean space, hence its projection minimizes length in the corresponding free homotopy class.  Its velocity is constant:
\[
  \operatorname{Length}(\gamma_P)
  =\int_0^1\norm{P}\dd t
  =\norm{P}
  =\sqrt{\inner{P}{P}}
  =\sqrt{\qheight(P)}.
\]
The stable-norm statement follows from homogeneity and continuity.
\end{proof}

\begin{corollary}[Primitive classes and multiple windings]\label{cor:primitive}
A nonzero $P\in M_E$ is primitive if and only if $\gamma_P$ is not an $m$-fold traversal of another closed geodesic for any $m>1$.  If $P=mQ$, then
\[
  \gamma_P(t)=\gamma_Q(mt\bmod1),
  \qquad
  \operatorname{Length}(\gamma_P)=m\operatorname{Length}(\gamma_Q),
\]
and $\qheight(mQ)=m^2\qheight(Q)$.
\end{corollary}

\subsection{Finite reconstruction of the height metric}

The topology alone records only $r$.  The marked length spectrum records the full height pairing.

\begin{theorem}[Marked-length reconstruction]\label{thm:marked-reconstruction}
Let $P_1,\ldots,P_r$ be a basis of $M_E$.  The $r(r+1)/2$ marked squared lengths
\[
  \operatorname{Length}(\gamma_{P_i})^2
  \quad(1\le i\le r),
  \qquad
  \operatorname{Length}(\gamma_{P_i+P_j})^2
  \quad(1\le i<j\le r)
\]
determine the complete Gram matrix $G_E(\mathcal B)$.  Explicitly,
\begin{equation}\label{eq:metric-reconstruction}
  \inner{P_i}{P_i}=\operatorname{Length}(\gamma_{P_i})^2,
\end{equation}
\begin{equation}\label{eq:cross-reconstruction}
  \inner{P_i}{P_j}
  =\frac12\left(
  \operatorname{Length}(\gamma_{P_i+P_j})^2
  -\operatorname{Length}(\gamma_{P_i})^2
  -\operatorname{Length}(\gamma_{P_j})^2
  \right).
\end{equation}
Consequently, the marked length data determine $\Reg(E/\Q)$ and the metric isometry class of $\T_E$.
\end{theorem}

\begin{proof}
Equation \eqref{eq:metric-reconstruction} is \eqref{eq:length-height}.  Equation \eqref{eq:cross-reconstruction} is the polarization identity \eqref{eq:height-pairing}.  A positive-definite Gram matrix determines the Euclidean lattice metric in the marked basis and therefore the quotient flat torus.
\end{proof}

\begin{remark}[Unmarked versus marked data]\label{rem:unmarked-spectrum}
The unmarked set of lengths need not identify which class produced which length, and spectral data need not determine a high-dimensional flat torus up to isometry.  Isospectral non-isometric flat tori exist \cite{Milnor1964,ConwaySloane1992,Schiemann1997,NilssonRowlettRydell2021}.  The arithmetic marking by $M_E$ is therefore substantive: it retains the group class attached to each geodesic.
\end{remark}

\section{Regulator, systoles, and finite-index coverings}

\subsection{Regulator as volume}

\begin{theorem}[Regulator--volume theorem]\label{thm:reg-volume}
The Riemannian volume of the height torus satisfies
\begin{equation}\label{eq:volume-regulator}
  \Vol(\T_E)=\sqrt{\Reg(E/\Q)},
  \qquad
  \Vol(\T_E)^2=\Reg(E/\Q).
\end{equation}
\end{theorem}

\begin{proof}
Let $P_1,\ldots,P_r$ be a lattice basis.  The volume of its fundamental parallelotope is $\sqrt{\det G_E}$.  This is the volume of the quotient $V_E/M_E$.  The rank-zero convention gives both sides equal to one.
\end{proof}

\begin{figure}[H]
  \centering
  \includegraphics[width=0.94\textwidth]{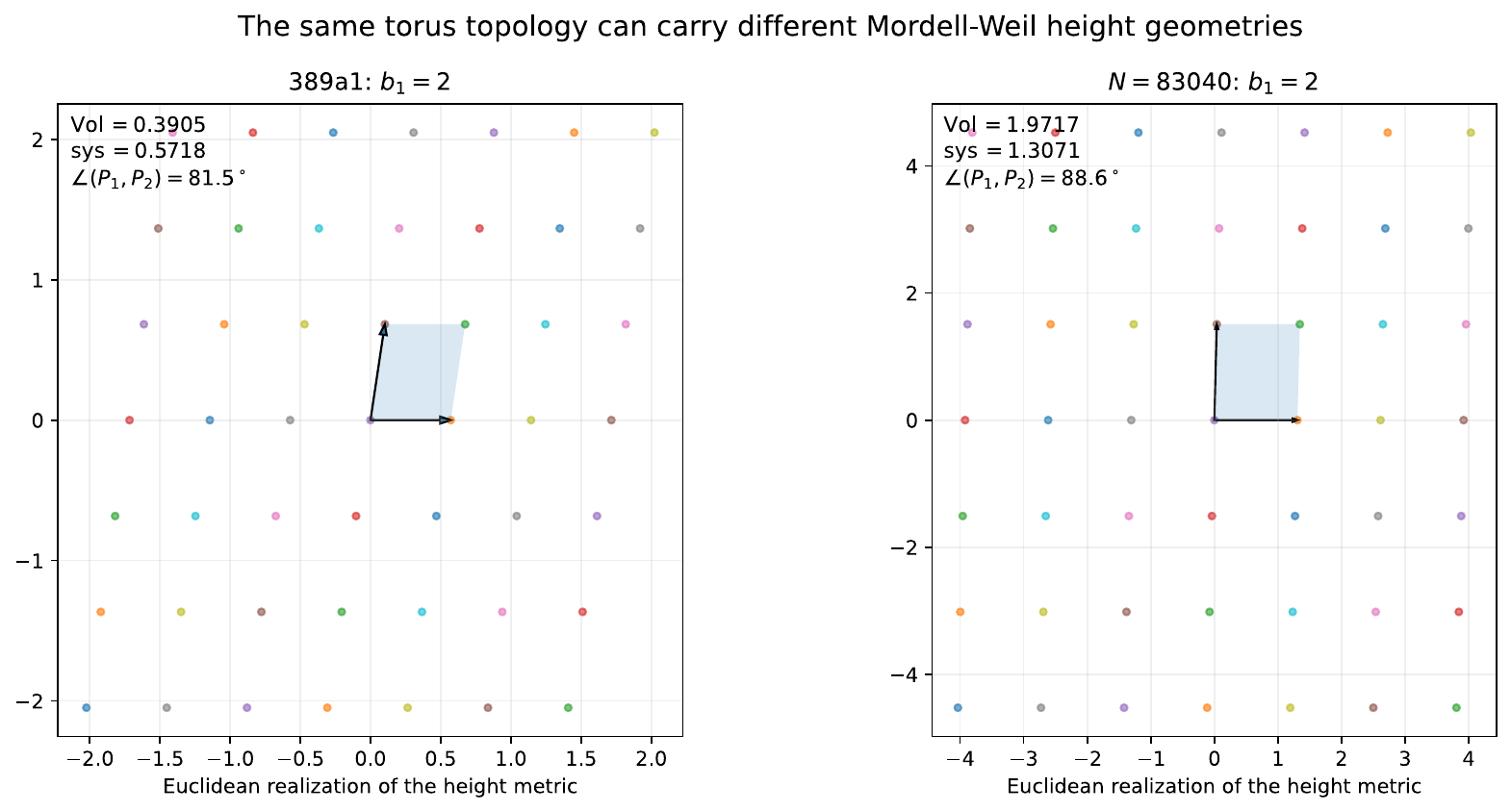}
  \caption{Two rank-two height lattices.  Their quotient tori have the same topology, $b_1=2$, but different N\'eron--Tate geometries, regulators, areas, and systoles.  The displayed fundamental parallelograms are Euclidean realizations of the computed Gram matrices.}
  \label{fig:rank2-lattices}
\end{figure}

\subsection{Systoles and successive minima}

Assume $r>0$.  Let $\lambda_1\le\cdots\le\lambda_r$ be the successive minima of $M_E$ with respect to the norm $\norm{P}^2=\qheight(P)$.

\begin{theorem}[Systole and minimal canonical height]\label{thm:systole}
The systole of $\T_E$ is
\begin{equation}\label{eq:systole-height}
  \sys(\T_E)^2
  =\min_{0\ne P\in M_E}\qheight(P)
  =\lambda_1^2.
\end{equation}
More generally, $\lambda_i^2$ is the least $H$ for which $M_E$ contains $i$ linearly independent classes of height at most $H$.
\end{theorem}

\begin{proof}
Every noncontractible loop represents a nonzero class $P\in H_1(\T_E,\Z)\cong M_E$.  By \cref{thm:height-length}, the shortest loop in that class has length $\sqrt{\qheight(P)}$.  Minimizing over nonzero classes proves the first statement.  The rest is the definition of successive minima translated through the same identity.  The terminology is compatible with the wider theory of systolic and filling inequalities \cite{Gromov1983}.
\end{proof}

Let $\kappa_r=\pi^{r/2}/\Gamma(r/2+1)$ be the volume of the unit ball in $\R^r$.

\begin{corollary}[Minkowski bounds in height form]\label{cor:minkowski}
The successive minima and the regulator satisfy
\begin{equation}\label{eq:minkowski-product}
  \frac{2^r}{r!\,\kappa_r}\sqrt{\Reg(E/\Q)}
  \le
  \prod_{i=1}^r\lambda_i
  \le
  \frac{2^r}{\kappa_r}\sqrt{\Reg(E/\Q)}.
\end{equation}
Equivalently,
\begin{equation}\label{eq:minkowski-regulator}
  \left(\frac{\kappa_r}{2^r}\right)^2
  \prod_{i=1}^r\lambda_i^2
  \le\Reg(E/\Q)\le
  \left(\frac{r!\,\kappa_r}{2^r}\right)^2
  \prod_{i=1}^r\lambda_i^2.
\end{equation}
\end{corollary}

\begin{proof}
Apply Minkowski's second theorem to the lattice $M_E$ and the unit ball of the height norm \cite{Minkowski1910,Cassels1959,ConwaySloane1999}.  Its covolume is $\sqrt{\Reg(E/\Q)}$ by \cref{thm:reg-volume}.  Related transference estimates between a lattice and its metric dual can be sharpened by results such as \cite{Banaszczyk1993}; they are relevant when comparing the shortest arithmetic geodesics with the lowest nonzero dual-lattice eigenvalues.
\end{proof}

\subsection{Saturation as a covering problem}

Computational point searches often produce a full-rank subgroup $L\subseteq M_E$ without immediately proving that $L=M_E$.  The index has a direct geometric interpretation.

\begin{theorem}[Finite-index covering theorem]\label{thm:index-cover}
Let $L\subseteq M_E$ be a full-rank sublattice of finite index
\[
  I=[M_E:L].
\]
Then the natural map
\begin{equation}\label{eq:index-cover-map}
  V_E/L\longrightarrow V_E/M_E=\T_E
\end{equation}
is an $I$-sheeted Riemannian covering.  Moreover,
\begin{equation}\label{eq:index-volume}
  \Vol(V_E/L)=I\Vol(\T_E),
\end{equation}
\begin{equation}\label{eq:index-regulator}
  \det G_L=I^2\Reg(E/\Q).
\end{equation}
\end{theorem}

\begin{proof}
Since $L\subseteq M_E$, quotienting by $L$ and then by the finite group $M_E/L$ gives $V_E/M_E$.  The deck group is $M_E/L$, of order $I$.  An $I$-sheeted local isometry multiplies total volume by $I$.  Squaring the volume identity and applying \cref{thm:reg-volume} gives \eqref{eq:index-regulator}.
\end{proof}

\begin{corollary}[Index detection from regulators]\label{cor:index-detection}
If the determinant of a candidate height Gram matrix is $R_L$ and the true regulator is independently known to be $R_E$, then
\[
  [M_E:L]=\sqrt{R_L/R_E}.
\]
Thus an unsaturated basis cannot be treated as a Mordell--Weil basis merely because its points are independent.
\end{corollary}

\section{Toroidal helices and four-dimensional windows}\label{sec:helices}

The term \emph{helix} is used here for a coordinate representation of a closed geodesic on a product of circles.  This preserves the valid geometric intuition while keeping the underlying homology class explicit.

\subsection{A \texorpdfstring{$2r$}{2r}-dimensional toroidal representation}

Choose a basis $\mathcal B=(P_1,\ldots,P_r)$ of $M_E$ and radii $R_i>0$.  Define
\begin{equation}\label{eq:toroidal-embedding}
  J_{\mathcal B,R}:\T_E\longrightarrow\R^{2r}
\end{equation}
by
\begin{align}\label{eq:toroidal-embedding-coordinates}
  J_{\mathcal B,R}(x_1,\ldots,x_r)
  :=(&R_1\cos2\pi x_1,R_1\sin2\pi x_1,\ldots,\\
     &R_r\cos2\pi x_r,R_r\sin2\pi x_r).
\end{align}
This is a smooth topological embedding of the underlying torus.  It is not generally an isometric embedding of the N\'eron--Tate metric because the height Gram matrix need not be diagonal in the selected basis.

\begin{proposition}[Toroidal-helical representative]\label{prop:helix}
If
\[
  P=n_1P_1+\cdots+n_rP_r,
\]
then
\begin{align}\label{eq:helix-formula}
  J_{\mathcal B,R}(\gamma_P(t))
  =(&R_1\cos2\pi n_1t,R_1\sin2\pi n_1t,\ldots,\\
     &R_r\cos2\pi n_rt,R_r\sin2\pi n_rt).
\end{align}
The winding vector $(n_1,\ldots,n_r)$ is the coordinate vector of the homology class $P$.  Under a change of Mordell--Weil basis, the winding vector transforms by $\mathrm{GL}_r(\Z)$.
\end{proposition}

\begin{proof}
In basis coordinates, $\gamma_P(t)=(n_1t,\ldots,n_rt)\bmod\Z^r$.  Substituting in \eqref{eq:toroidal-embedding-coordinates} gives the formula.  The transformation law is the ordinary integral change of lattice coordinates.
\end{proof}

\begin{figure}[H]
  \centering
  \includegraphics[width=0.68\textwidth]{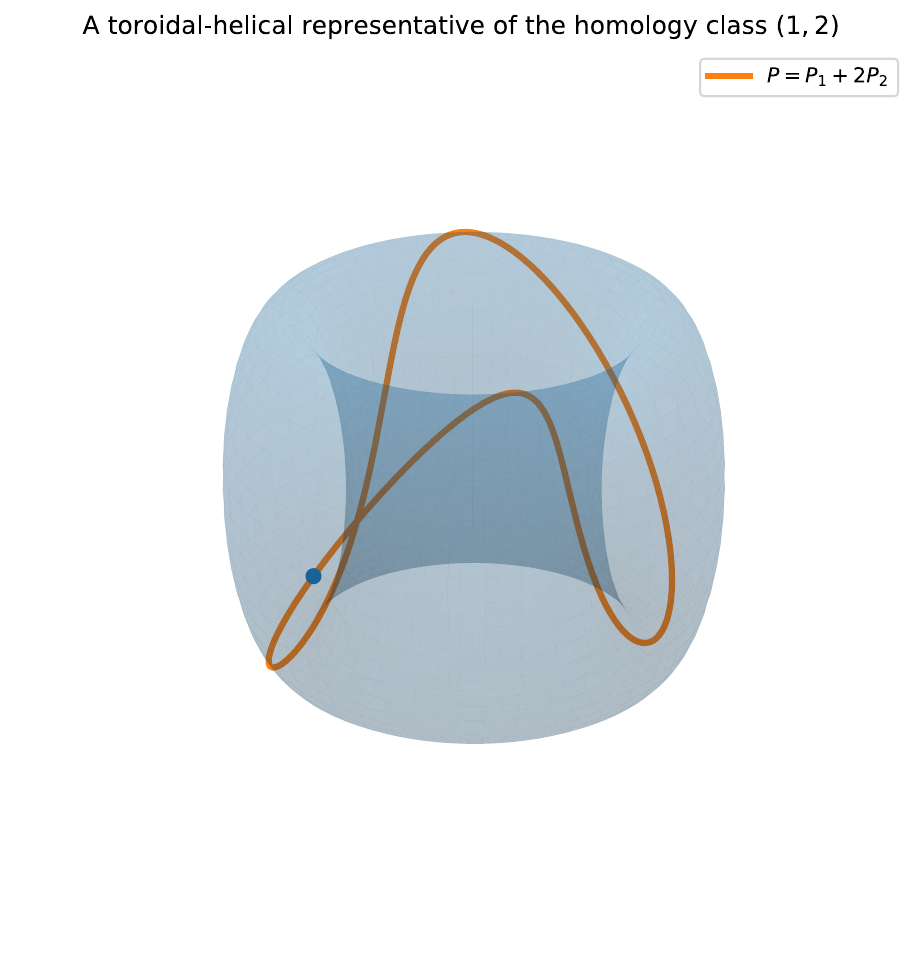}
  \caption{A toroidal-helical representative of the class $P=P_1+2P_2$ on a rank-two torus.  The drawing represents the winding data $(1,2)$; its Euclidean arc length is not asserted to equal the canonical height unless the ambient embedding is chosen isometrically.}
  \label{fig:helix}
\end{figure}

\subsection{Character windows in four dimensions}

Let
\[
  M_E^\vee:=\Hom(M_E,\Z)
\]
be the integral character lattice.  For $\ell_1,\ell_2\in M_E^\vee$, define
\begin{equation}\label{eq:4d-window}
  W_{\ell_1,\ell_2}:\T_E\longrightarrow(S^1)^2\subset\R^4
\end{equation}
by
\begin{align}\label{eq:4d-window-formula}
  W_{\ell_1,\ell_2}(x)
  :=(&\cos2\pi\ell_1(x),\sin2\pi\ell_1(x),\\
     &\cos2\pi\ell_2(x),\sin2\pi\ell_2(x)).
\end{align}

\begin{theorem}[Four-dimensional window theorem]\label{thm:4d-window}
The induced map on first homology is
\begin{equation}\label{eq:4d-homology-map}
  (W_{\ell_1,\ell_2})_*:M_E\longrightarrow\Z^2,
  \qquad
  P\longmapsto(\ell_1(P),\ell_2(P)).
\end{equation}
A single four-dimensional character window therefore retains at most two independent homology directions.  If $r=2$ and $\ell_1,\ell_2$ form the dual basis, the window gives the standard embedding of $\T_E$ in $\R^4$.  If $r>2$, the homological kernel has rank at least $r-2$.
\end{theorem}

\begin{proof}
A loop $\gamma_P$ winds $\ell_i(P)$ times around the $i$th target circle.  This gives \eqref{eq:4d-homology-map}.  A homomorphism $\Z^r\to\Z^2$ has kernel rank at least $r-2$.
\end{proof}

\begin{theorem}[Minimum number of complete four-dimensional windows]\label{thm:min-windows}
Within the character-window family, the minimum number of four-dimensional windows whose combined homology map is injective is
\begin{equation}\label{eq:min-windows}
  m_{\min}=\left\lceil\frac r2\right\rceil.
\end{equation}
\end{theorem}

\begin{proof}
A collection of $m$ windows induces a homomorphism $\Z^r\to\Z^{2m}$.  Injectivity requires $r\le2m$, so $m\ge\lceil r/2\rceil$.  Conversely, choose a dual basis $\ell_1,\ldots,\ell_r$ and distribute the characters in pairs among $\lceil r/2\rceil$ windows, taking the final unused character to be zero if $r$ is odd.  The combined map records all $r$ coordinates and is injective.
\end{proof}

\begin{figure}[H]
  \centering
  \includegraphics[width=0.95\textwidth]{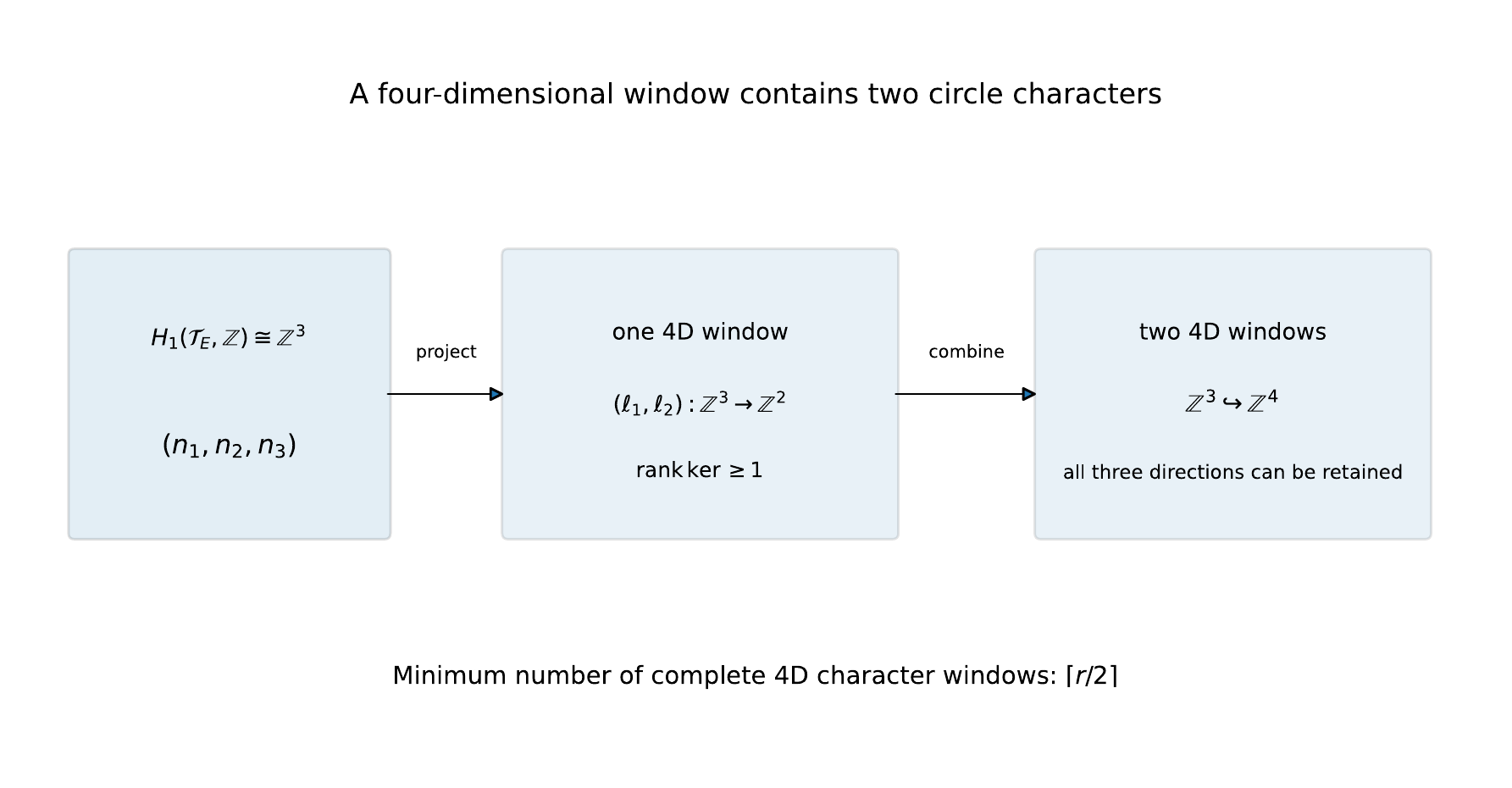}
  \caption{A four-dimensional window consists of two circle characters.  For rank three, one window necessarily loses at least one homology direction, whereas two windows can retain all three.}
  \label{fig:4d-windows}
\end{figure}

\begin{remark}
The theorem is about this explicit character-window construction, not the unrestricted minimal Euclidean embedding dimension of an abstract torus.  It gives a precise interpretation of what a four-dimensional visualization can and cannot claim about arithmetic rank.
\end{remark}

\section{Theta series and spectral geometry}\label{sec:spectral}

The height lattice has two closely related analytic objects: its theta series and the heat trace of its quotient torus.  These are rigorous and computable; they must not be confused with the Hasse--Weil $L$-function.

\subsection{Height theta series}

Define
\begin{equation}\label{eq:theta}
  \Theta_E(u):=\sum_{P\in M_E}\exp\bigl(-\pi u\qheight(P)\bigr),
  \qquad u>0.
\end{equation}
Let the metric dual lattice be
\begin{equation}\label{eq:metric-dual}
  M_E^*:=\{\xi\in V_E:\inner{\xi}{P}\in\Z\ \text{for all }P\in M_E\}.
\end{equation}

\begin{theorem}[Poisson identity]\label{thm:poisson}
The theta series satisfies
\begin{equation}\label{eq:poisson}
  \Theta_E(u)
  =\frac{u^{-r/2}}{\Vol(\T_E)}
  \sum_{\xi\in M_E^*}
  \exp\left(-\frac{\pi\norm{\xi}^2}{u}\right).
\end{equation}
\end{theorem}

\begin{proof}
Apply Poisson summation to the Gaussian $x\mapsto\exp(-\pi u\norm{x}^2)$ on $V_E$.  The Fourier transform is $u^{-r/2}\exp(-\pi\norm{\xi}^2/u)$, and the lattice covolume is $\Vol(\T_E)$.
\end{proof}

The theta series counts Mordell--Weil classes by canonical height, with $P$ and $-P$ counted separately unless $P=0$.  It therefore packages more metric information than rank and regulator alone.  This lies in the broader spectral-geometric question of how much geometry can be recovered from spectral data, classically emphasized by Kac \cite{Kac1966}.  Nonetheless, theta series can fail to distinguish nonisometric lattices in sufficiently high dimension \cite{ConwaySloane1992,Schiemann1997}; it is an invariant, not a complete invariant in all ranks.

\begin{figure}[H]
  \centering
  \includegraphics[width=0.82\textwidth]{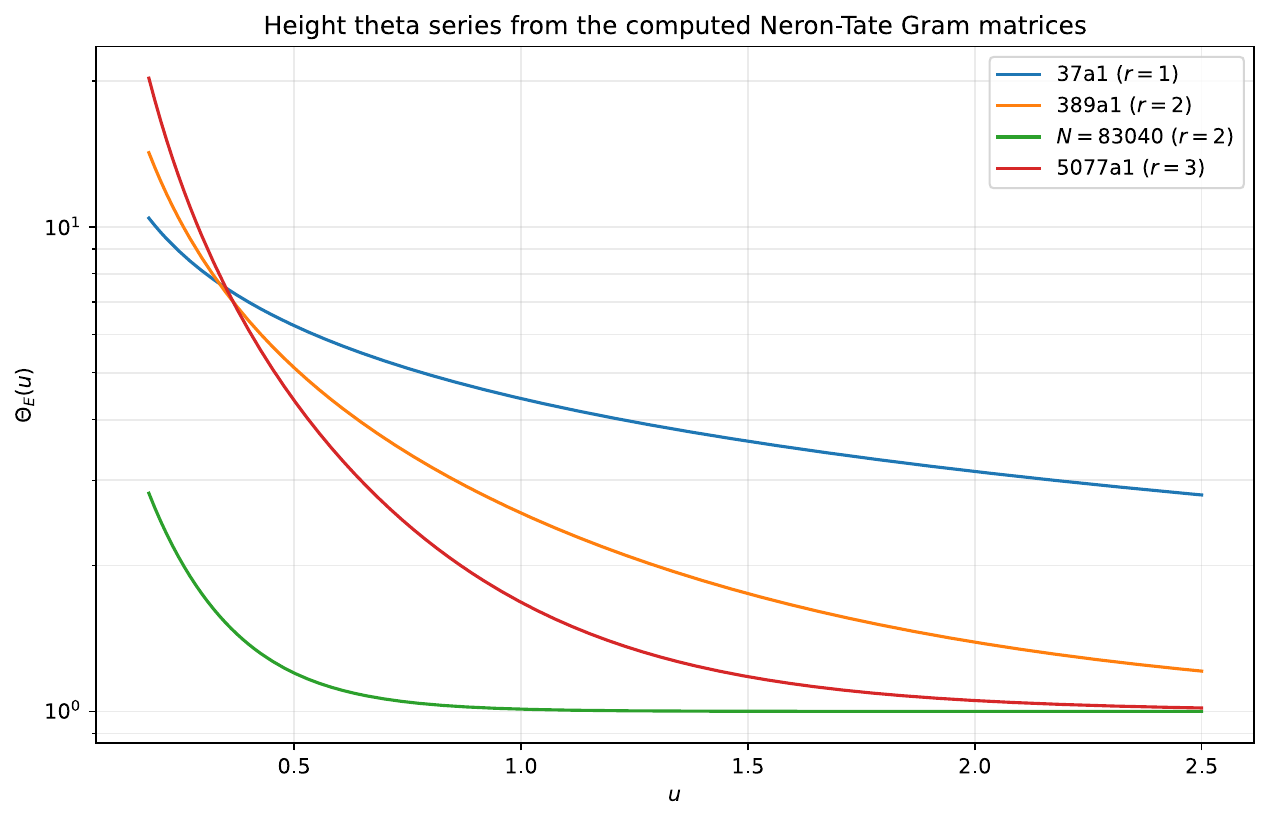}
  \caption{Finite truncations of the height theta series for four computed Mordell--Weil lattices.  The curves have different ranks and height geometries; the plotted sums use the Gram matrices listed in \cref{sec:computations}.}
  \label{fig:theta}
\end{figure}

\subsection{Laplacian and heat trace}

Let $\Delta_E$ be the nonnegative scalar Laplace--Beltrami operator on $\T_E$.

\begin{theorem}[Flat-torus spectrum]\label{thm:flat-spectrum}
The functions
\[
  e_\xi(x):=\exp\bigl(2\pi i\inner{\xi}{x}\bigr),
  \qquad \xi\in M_E^*,
\]
form an orthogonal eigenbasis and
\begin{equation}\label{eq:spectrum}
  \Delta_Ee_\xi=4\pi^2\norm{\xi}^2e_\xi.
\end{equation}
Hence
\begin{equation}\label{eq:heat-dual}
  \Tr(e^{-t\Delta_E})
  =\sum_{\xi\in M_E^*}\exp(-4\pi^2t\norm{\xi}^2).
\end{equation}
\end{theorem}

\begin{proof}
The exponentials are precisely the characters of the compact abelian group $V_E/M_E$.  Differentiation in orthonormal coordinates gives the eigenvalue $4\pi^2\norm{\xi}^2$.  Fourier analysis on a flat torus yields completeness \cite{BergerGauduchonMazet1971,Chavel1984}.
\end{proof}

\begin{corollary}[Exact heat-trace formula]\label{cor:heat-trace}
For every $t>0$,
\begin{equation}\label{eq:heat-primal}
  \Tr(e^{-t\Delta_E})
  =\frac{\Vol(\T_E)}{(4\pi t)^{r/2}}
  \sum_{P\in M_E}\exp\left(-\frac{\qheight(P)}{4t}\right).
\end{equation}
\end{corollary}

\begin{proof}
Apply \cref{thm:poisson} with $u=4\pi t$ after matching the Gaussian normalization, or apply Poisson summation directly to \eqref{eq:heat-dual}.
\end{proof}

\begin{corollary}[Recovery of rank, regulator, and systole]\label{cor:spectral-recovery}
As $t\downarrow0$,
\begin{equation}\label{eq:heat-asymptotic}
  \Tr(e^{-t\Delta_E})
  =\frac{\sqrt{\Reg(E/\Q)}}{(4\pi t)^{r/2}}
  \left(1+O\left(\exp\left[-\frac{\sys(\T_E)^2}{4t}\right]\right)\right).
\end{equation}
Consequently,
\begin{equation}\label{eq:recover-rank}
  r=-2\lim_{t\downarrow0}
  \frac{\log\Tr(e^{-t\Delta_E})}{\log t},
\end{equation}
\begin{equation}\label{eq:recover-volume}
  \sqrt{\Reg(E/\Q)}
  =\lim_{t\downarrow0}(4\pi t)^{r/2}\Tr(e^{-t\Delta_E}).
\end{equation}
If $r>0$, then
\begin{equation}\label{eq:recover-systole}
  \sys(\T_E)^2
  =\lim_{t\downarrow0}-4t\log\left(
  \frac{(4\pi t)^{r/2}}{\Vol(\T_E)}\Tr(e^{-t\Delta_E})-1
  \right).
\end{equation}
\end{corollary}

\begin{proof}
In \eqref{eq:heat-primal}, the zero vector contributes one.  The smallest nonzero exponent is $\sys(\T_E)^2/(4t)$ by \cref{thm:systole}.  The stated limits follow by taking logarithms and leading terms.
\end{proof}

\begin{figure}[H]
  \centering
  \includegraphics[width=0.82\textwidth]{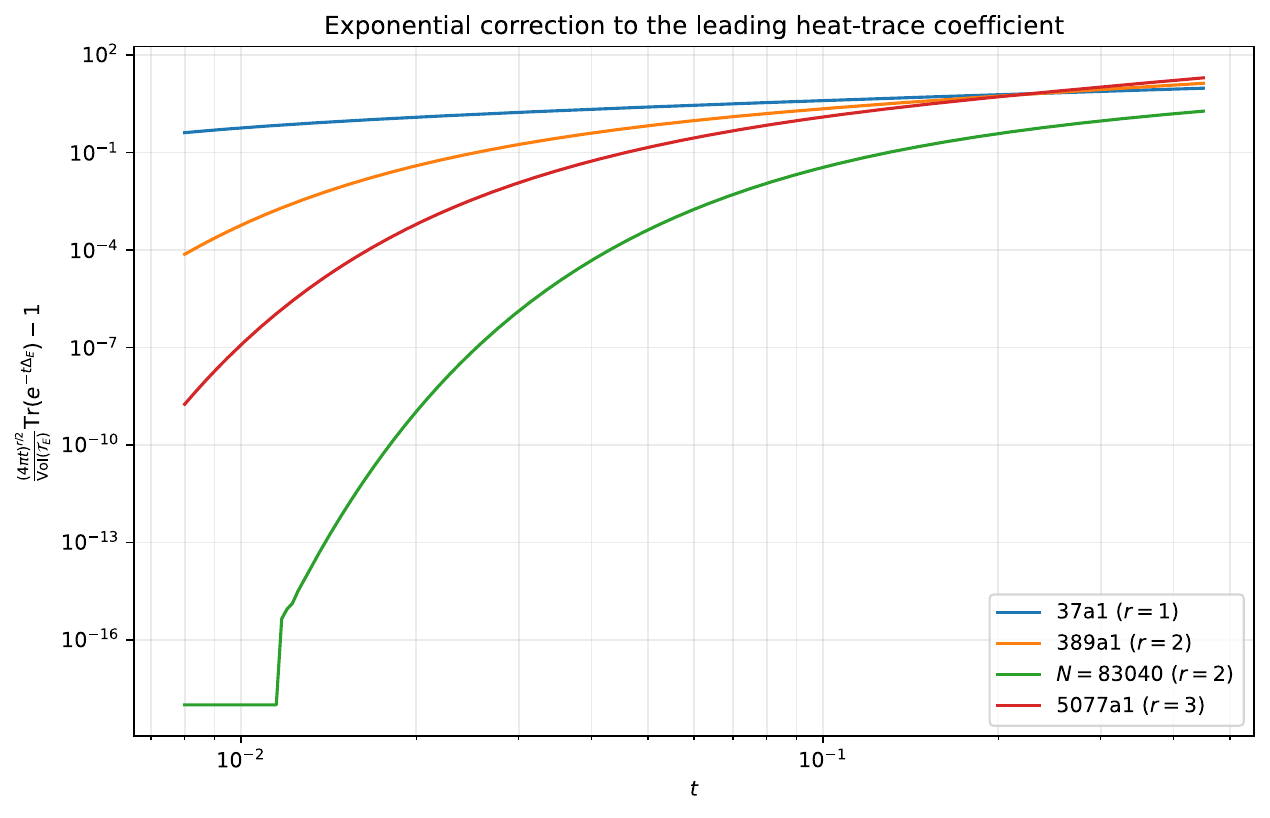}
  \caption{Finite-lattice evaluation of the normalized heat-trace correction
  $(4\pi t)^{r/2}\Tr(e^{-t\Delta_E})/\Vol(\T_E)-1$.  The small-$t$ slope is controlled by the shortest nonzero canonical height.}
  \label{fig:heat}
\end{figure}

\subsection{Hodge-theoretic zero modes}

Let $\Delta_E^{(k)}$ denote the Hodge Laplacian on $k$-forms.  Hodge theory gives
\begin{equation}\label{eq:Hodge}
  \ker\Delta_E^{(k)}\cong H^k(\T_E,\R).
\end{equation}
Thus
\begin{equation}\label{eq:harmonic-one-forms}
  \dim\ker\Delta_E^{(1)}=b_1(\T_E)=r.
\end{equation}
This equality is intrinsic and rigorous.  The open issue is not the zero-mode count of $\T_E$; it is the independent analytic identification of that count with the central multiplicity of $L(E,s)$.

\section{Behavior under isogenies}

Let $\phi:E\to E'$ be an isogeny of degree $d$.  It induces an injective homomorphism on free Mordell--Weil groups with finite cokernel.  Write
\[
  I_\phi:=[M_{E'}:\phi(M_E)].
\]

\begin{theorem}[Isogeny scaling]\label{thm:isogeny}
For $P,Q\in M_E$,
\begin{equation}\label{eq:isogeny-height}
  q_{E'}(\phi P)=d\,q_E(P),
  \qquad
  \langle\phi P,\phi Q\rangle_{E'}=d\langle P,Q\rangle_E.
\end{equation}
Consequently,
\begin{equation}\label{eq:isogeny-length}
  \operatorname{Length}_{\T_{E'}}(\gamma_{\phi P})
  =\sqrt d\,\operatorname{Length}_{\T_E}(\gamma_P),
\end{equation}
\begin{equation}\label{eq:isogeny-reg}
  d^r\Reg(E/\Q)=I_\phi^2\Reg(E'/\Q),
\end{equation}
\begin{equation}\label{eq:isogeny-volume}
  d^{r/2}\Vol(\T_E)=I_\phi\Vol(\T_{E'}).
\end{equation}
\end{theorem}

\begin{proof}
The canonical height associated with the origin divisor scales by the degree of an isogeny \cite{Silverman2009,Silverman1994}.  This proves \eqref{eq:isogeny-height} and \eqref{eq:isogeny-length}.  If $P_1,\ldots,P_r$ is a basis of $M_E$, the Gram determinant of $\phi(M_E)$ in $M_{E'}$ equals $d^r\Reg(E/\Q)$.  By the finite-index formula \eqref{eq:index-regulator}, it also equals $I_\phi^2\Reg(E'/\Q)$.  Taking square roots gives the volume identity.
\end{proof}

The theorem shows that the height-torus construction behaves predictably under one of the central equivalence relations in the arithmetic of elliptic curves.  It is not invariant under isogeny as a metric torus, nor should it be; its scaling and finite-index defect are exactly controlled.

\section{An independent analytic side of the elliptic-curve \texorpdfstring{$L$}{L}-function}\label{sec:analytic-side}

The height torus gives a rigorous metric-topological object, but it is built from the global Mordell--Weil group.  To formulate a non-circular comparison, the analytic object must be constructed independently from the local factors, conductor, and Archimedean completion of $L(E,s)$.  This section fixes that second construction.

\subsection{Completion, local logarithmic coefficients, and convention}

We use the Fourier transform
\begin{equation}\label{eq:residual-fourier}
  \widehat h(\xi)=\int_{\R}h(u)\e^{-2\pi i u\xi}\dd u,
  \qquad
  h(u)=\int_{\R}\widehat h(\xi)\e^{2\pi i u\xi}\dd\xi.
\end{equation}
The completed elliptic-curve $L$-function is
\begin{equation}\label{eq:completed-E-new}
  \Lambda(E,s)=N_E^{s/2}(2\pi)^{-s}\Gamma(s)L(E,s),
\end{equation}
with functional equation $\Lambda(E,s)=w_E\Lambda(E,2-s)$.  For every prime write the local factor as
\begin{equation}\label{eq:local-roots-general}
  L_p(E,s)=\prod_{j=1}^{d_p}(1-\alpha_{p,j}p^{-s})^{-1},
  \qquad d_p\in\{0,1,2\},
\end{equation}
where $d_p=2$ at good primes.  Define the generalized von Mangoldt coefficients by
\begin{equation}\label{eq:Lambda-E}
  \Lambda_E(p^m)=(\log p)\sum_{j=1}^{d_p}\alpha_{p,j}^{m},
  \qquad m\ge1,
\end{equation}
and $\Lambda_E(n)=0$ when $n$ is not a prime power.  Then, in the half-plane of absolute convergence,
\begin{equation}\label{eq:log-derivative-E}
  -\frac{L'}{L}(E,s)=\sum_{n\ge1}\frac{\Lambda_E(n)}{n^s}.
\end{equation}
At a good prime, $\alpha_p+\beta_p=a_p$ and $\alpha_p\beta_p=p$, so the quantities $\alpha_p^m+\beta_p^m$ can be computed recursively from $a_p$ and $p$ without choosing complex square roots.

\subsection{A centered explicit-formula functional}

Let $\mathscr H$ denote the real vector space of even entire functions $h$ satisfying $h(\overline z)=\overline{h(z)}$ and, for every $A,Y>0$,
\begin{equation}\label{eq:explicit-test-class}
  \sup_{\substack{x\in\R\\ |y|\le Y}}
  (1+|x|)^A|h(x+iy)|<\infty.
\end{equation}
The zero sum is taken symmetrically and zeros are counted with multiplicity.

\begin{theorem}[Centered explicit formula]\label{thm:elliptic-explicit-formula}
For $h\in\mathscr H$ whose Fourier transform is compactly supported,
\begin{align}\label{eq:elliptic-explicit-formula}
  \sum_{\rho}h\!\left(\frac{\rho-1}{i}\right)
  ={}&\frac{1}{2\pi}\int_{\R}h(u)
  \left[
    \log N_E-2\log(2\pi)
    +2\operatorname{Re}\frac{\Gamma'}{\Gamma}(1+iu)
  \right]\dd u \\
  &-\frac{1}{\pi}\sum_{n\ge1}\frac{\Lambda_E(n)}{n}
  \widehat h\!\left(\frac{\log n}{2\pi}\right).
\end{align}
Here $\rho$ ranges over the nontrivial zeros of $L(E,s)$.  The right-hand side depends only on the conductor, the Archimedean factor, and the local Euler data of $E$.
\end{theorem}

\begin{proof}
Choose $c>3/2$ and integrate the logarithmic derivative of the completed function against $h((s-1)/i)$ around a rectangle with vertical sides $\operatorname{Re}s=c$ and $\operatorname{Re}s=2-c$.  The completed function is entire of order one, so the residue theorem gives the symmetric zero sum.  The decay of $h$ makes the horizontal integrals tend to zero.  On the left vertical line, use the functional equation to reflect the logarithmic derivative to the right line.  Substituting
\[
  \frac{\Lambda'}{\Lambda}(E,s)
  =\frac12\log N_E-\log(2\pi)
  +\frac{\Gamma'}{\Gamma}(s)+\frac{L'}{L}(E,s)
\]
produces the conductor and Gamma integral.  On the right line, expand $-L'/L$ by \eqref{eq:log-derivative-E}; normal convergence permits termwise integration.  Fourier inversion \eqref{eq:residual-fourier} evaluates the $n$th integral at $\log n/(2\pi)$ and the two reflected vertical lines combine to give the factor $1/\pi$.  This is the standard Weil explicit-formula contour argument specialized to the degree-two completion \eqref{eq:completed-E-new}; see \cite{Weil1952,IwaniecKowalski2004} for the general framework.
\end{proof}

\begin{definition}[Explicit-formula analytic functional]\label{def:canonical-analytic-functional}
For admissible $h$, let $\mathcal E_E(h)$ denote either side of \eqref{eq:elliptic-explicit-formula}.  This functional is independently defined once the completion, Fourier convention, and test class are fixed.  Its decomposition into Archimedean and finite-place parts is written
\begin{equation}\label{eq:analytic-decomposition}
  \mathcal E_E(h)=\mathcal A_E(h)+\mathcal P_E(h).
\end{equation}
This functional is defined before any operator on the height torus is introduced.
\end{definition}

\subsection{Exact support-induced finiteness}

\begin{theorem}[Exact finite prime-power side]\label{thm:exact-prime-support}
If $\operatorname{supp}\widehat h\subseteq[-\Delta,\Delta]$, then the finite-place term $\mathcal P_E(h)$ contains only prime powers satisfying
\begin{equation}\label{eq:prime-support-bound}
  n\le \exp(2\pi\Delta).
\end{equation}
There is no omitted prime tail.
\end{theorem}

\begin{proof}
The $n$th local term in \eqref{eq:elliptic-explicit-formula} samples $\widehat h$ at $\log n/(2\pi)$.  This sample vanishes whenever $\log n/(2\pi)>\Delta$.  The remaining set of integers is finite.
\end{proof}

\begin{figure}[H]
  \centering
  \includegraphics[width=0.96\textwidth]{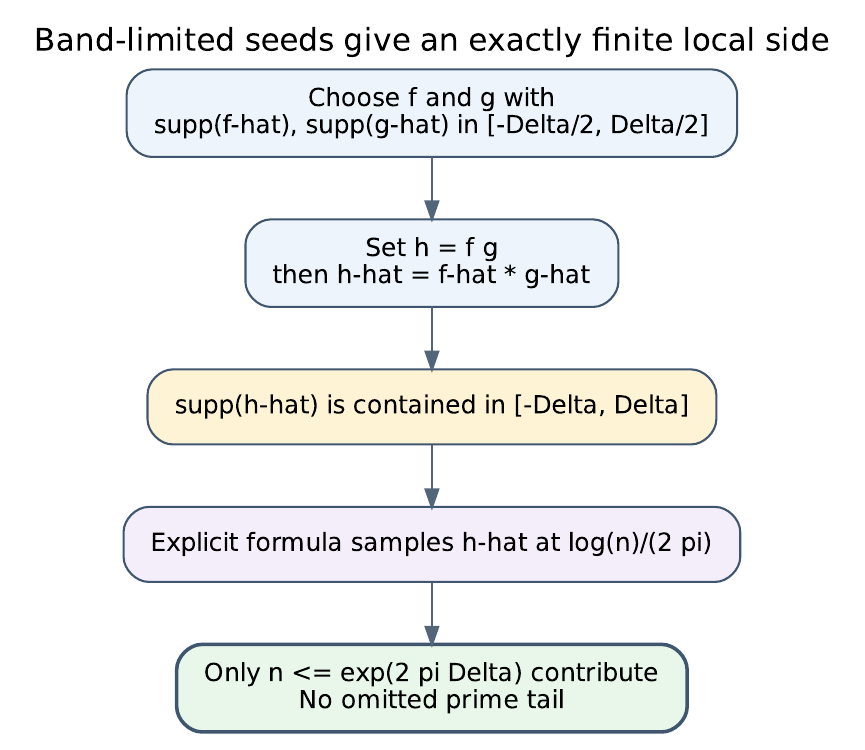}
  \caption{Compact Fourier support makes the prime-power side exactly finite.  The resulting computation still requires the full Archimedean term and a validated quadrature enclosure; support truncation alone is not a positivity theorem.}
  \label{fig:exact-support}
\end{figure}

\subsection{An analytic form with an exactly finite local side}

For $\Delta>0$, define the real even Paley--Wiener--Schwartz seed space
\begin{equation}\label{eq:test-seed-space}
  \mathscr P_{\Delta/2}
  :=\left\{f=\mathcal F^{-1}\varphi:
  \varphi\in C_{c,\mathrm{even}}^\infty
  \bigl((-\Delta/2,\Delta/2);\R\bigr)\right\}.
\end{equation}
\begin{lemma}[Regularity of the seed class]\label{lem:seed-regularity}
Every $f\in\mathscr P_{\Delta/2}$ is real and even on $\R$, extends to an entire function of exponential type, and is rapidly decreasing uniformly on each horizontal strip.  Hence $f,g\in\mathscr P_{\Delta/2}$ implies $fg\in\mathscr H$.
\end{lemma}

\begin{proof}
Write $f(z)=\int\widehat f(\xi)e^{2\pi i z\xi}\dd\xi$.  Compact support permits differentiation under the integral and gives an entire extension with exponential growth controlled by $|\operatorname{Im}z|$.  Repeated integration by parts in $\xi$ gives arbitrary inverse powers of $|\operatorname{Re}z|$, uniformly when $\operatorname{Im}z$ remains in a fixed compact interval.  Reality and evenness follow from the corresponding properties of $\widehat f$.  Products preserve these properties and the stripwise rapid decay.
\end{proof}

For $f,g\in\mathscr P_{\Delta/2}$ set
\begin{equation}\label{eq:seed-test-map}
  h_{f,g}(u):=f(u)g(u).
\end{equation}
Then
\begin{equation}\label{eq:seed-convolution}
  \widehat h_{f,g}=\widehat f*\widehat g,
  \qquad
  \operatorname{supp}\widehat h_{f,g}
  \subseteq[-\Delta,\Delta].
\end{equation}
We define the real symmetric bilinear form
\begin{equation}\label{eq:analytic-form}
  W^{\mathrm{an}}_{E,\Delta}(f,g)
  :=\mathcal E_E(h_{f,g}).
\end{equation}

\begin{proposition}[Independence, symmetry, and finite local assembly]\label{prop:analytic-form-properties}
The form $W^{\mathrm{an}}_{E,\Delta}$ is symmetric and is constructed from the completed $L$-function without using a Mordell--Weil basis, canonical heights, the regulator, or the algebraic rank.  For every finite family $f_1,\ldots,f_N\in\mathscr P_{\Delta/2}$, its matrix
\begin{equation}\label{eq:analytic-matrix}
  G^{\mathrm{an}}_{E,\Delta}
  =\bigl(W^{\mathrm{an}}_{E,\Delta}(f_i,f_j)\bigr)_{i,j=1}^N
\end{equation}
has an exactly finite prime-power contribution.
\end{proposition}

\begin{proof}
Symmetry follows from $h_{f,g}=h_{g,f}$ and the linearity of $\mathcal E_E$.  Equation \eqref{eq:seed-convolution} follows from the Fourier product--convolution identity, and the support inclusion is the Minkowski-sum bound for convolution supports.  The input data in \eqref{eq:elliptic-explicit-formula} are $N_E$, the Gamma factor, and $\Lambda_E(n)$; none is defined from $M_E$ or its height form.  Exact finiteness of the local contribution follows from \cref{thm:exact-prime-support}.
\end{proof}

\begin{warning}[Positivity is not built into the definition]\label{warn:analytic-positivity}
The form $W^{\mathrm{an}}_{E,\Delta}$ is independently defined relative to the fixed explicit formula and seed-to-test map.  It is not declared positive by definition.  Any sign theorem must be proved from the analytic zero distribution or by an independent comparison.
\end{warning}

\section{The height-torus residual and spectral obstructions}\label{sec:residual}

We now place the analytic construction beside a positive form generated by the height-torus spectrum.  Their difference, rather than the automatic positivity of the torus form, is the comparison object.

\subsection{The compact torus trace form}

Let
\begin{equation}\label{eq:D-one-form}
  D_E^{(1)}:=\sqrt{\Delta_E^{(1)}}
\end{equation}
be the nonnegative square root of the Hodge Laplacian on one-forms.  For a declared scale $\tau>0$ and $f,g\in\mathscr P_{\Delta/2}$, define
\begin{equation}\label{eq:torus-form}
  q^{\mathrm{tor}}_{E,\Delta,\tau}(f,g)
  :=\Tr\!\left[f(\tau D_E^{(1)})g(\tau D_E^{(1)})\right].
\end{equation}
Since $f$ and $g$ are Schwartz on $\R$ and the flat-torus counting function has polynomial growth, the product multiplier is trace class.  The dual-lattice description of the Hodge spectrum gives the absolutely convergent expansion
\begin{equation}\label{eq:torus-form-expanded}
  q^{\mathrm{tor}}_{E,\Delta,\tau}(f,g)
  =r\sum_{\xi\in M_E^*}
  f(2\pi\tau\norm{\xi})g(2\pi\tau\norm{\xi}),
\end{equation}
with the rank-zero case interpreted as the zero form.

\begin{proposition}[Trace class and positivity of the torus form]\label{prop:torus-form-psd}
The form $q^{\mathrm{tor}}_{E,\Delta,\tau}$ is well defined, symmetric, and positive semidefinite.  It depends on the height lattice and the explicit scale $\tau$, but not on the zeros or analytic rank of $L(E,s)$.
\end{proposition}

\begin{proof}
For every $A>0$, the Schwartz decay gives
\[
  \abs{f(x)g(x)}\le C_A(1+\abs{x})^{-A}.
\]
The number of dual-lattice vectors with norm at most $X$ is $O(X^r)$, so choosing $A>r$ proves absolute convergence of \eqref{eq:torus-form-expanded} and trace-classness of the product multiplier.  For $f=g$,
\[
  q^{\mathrm{tor}}_{E,\Delta,\tau}(f,f)
  =r\sum_{\xi\in M_E^*}f(2\pi\tau\norm{\xi})^2\ge0.
\]
The spectrum is determined by $G_E^{-1}$ and contains no analytic $L$-function input.
\end{proof}

\subsection{Definition and basis covariance of the residual}

\begin{definition}[Height-torus residual]\label{def:height-torus-residual}
For fixed $(E,\Delta,\tau)$ define
\begin{equation}\label{eq:residual-form}
  \mathcal R_{E,\Delta,\tau}
  :=W^{\mathrm{an}}_{E,\Delta}-q^{\mathrm{tor}}_{E,\Delta,\tau}.
\end{equation}
On a finite family $\mathcal F=(f_1,\ldots,f_N)$, write
\begin{equation}\label{eq:residual-matrix}
  R_{E,\Delta,\tau}(\mathcal F)
  =G^{\mathrm{an}}_{E,\Delta}(\mathcal F)
  -G^{\mathrm{tor}}_{E,\Delta,\tau}(\mathcal F).
\end{equation}
\end{definition}

\begin{proposition}[Central-mode decomposition]\label{prop:central-mode-decomposition}
Let $r_{\mathrm{an}}(E)=\ord_{s=1}L(E,s)$.  For $f,g\in\mathscr P_{\Delta/2}$,
\begin{equation}\label{eq:central-mode-decomposition}
  \mathcal R_{E,\Delta,\tau}(f,g)
  =\bigl(r_{\mathrm{an}}(E)-r\bigr)f(0)g(0)
   +\mathcal R^{\circ}_{E,\Delta,\tau}(f,g),
\end{equation}
where
\begin{align}\label{eq:noncentral-residual}
  \mathcal R^{\circ}_{E,\Delta,\tau}(f,g)
  :={}&\sum_{\substack{\rho\\ \rho\ne1}}
  f\!\left(\frac{\rho-1}{i}\right)
  g\!\left(\frac{\rho-1}{i}\right) \\
  &-r\sum_{\substack{\xi\in M_E^*\\ \xi\ne0}}
  f(2\pi\tau\norm{\xi})g(2\pi\tau\norm{\xi}).
\end{align}
Thus the explicitly separated zero-frequency contribution to the residual is exactly the analytic-rank discrepancy multiplied by $(f,g)\mapsto f(0)g(0)$.
\end{proposition}

\begin{proof}
A zero $\rho=1$ contributes $h_{f,g}(0)=f(0)g(0)$ to the zero side of \eqref{eq:elliptic-explicit-formula}; its multiplicity is $r_{\mathrm{an}}(E)$.  In \eqref{eq:torus-form-expanded}, the dual-lattice vector $\xi=0$ contributes $r f(0)g(0)$ because the harmonic one-form space has dimension $r$.  Separating these two zero-frequency terms from the remaining sums gives \eqref{eq:central-mode-decomposition}.
\end{proof}

\begin{remark}[What the decomposition proves]\label{rem:central-decomposition-scope}
Equation \eqref{eq:central-mode-decomposition} exposes the rank discrepancy inside the residual without assuming BSD.  It does not show that the noncentral residual vanishes, has a sign, or is independent of $\tau$.  Any such statement requires an additional theorem or a validated finite computation with tail control.
\end{remark}

\begin{theorem}[Congruence covariance]\label{thm:residual-covariance}
Let $A\in\operatorname{GL}_N(\R)$ and define a new test family by $g_i=\sum_j f_jA_{ji}$.  Then
\begin{equation}\label{eq:residual-congruence}
  R_{E,\Delta,\tau}(\mathcal G)
  =A^{\mathsf T}R_{E,\Delta,\tau}(\mathcal F)A.
\end{equation}
Consequently, exact vanishing, rank, and inertia are independent of the chosen basis of the same finite test space.
\end{theorem}

\begin{proof}
Both $W^{\mathrm{an}}$ and $q^{\mathrm{tor}}$ are bilinear.  Their matrices therefore transform by congruence under a basis change.  Subtraction proves \eqref{eq:residual-congruence}, and Sylvester's law of inertia gives the final statement.
\end{proof}

\begin{figure}[H]
  \centering
  \includegraphics[width=0.99\textwidth]{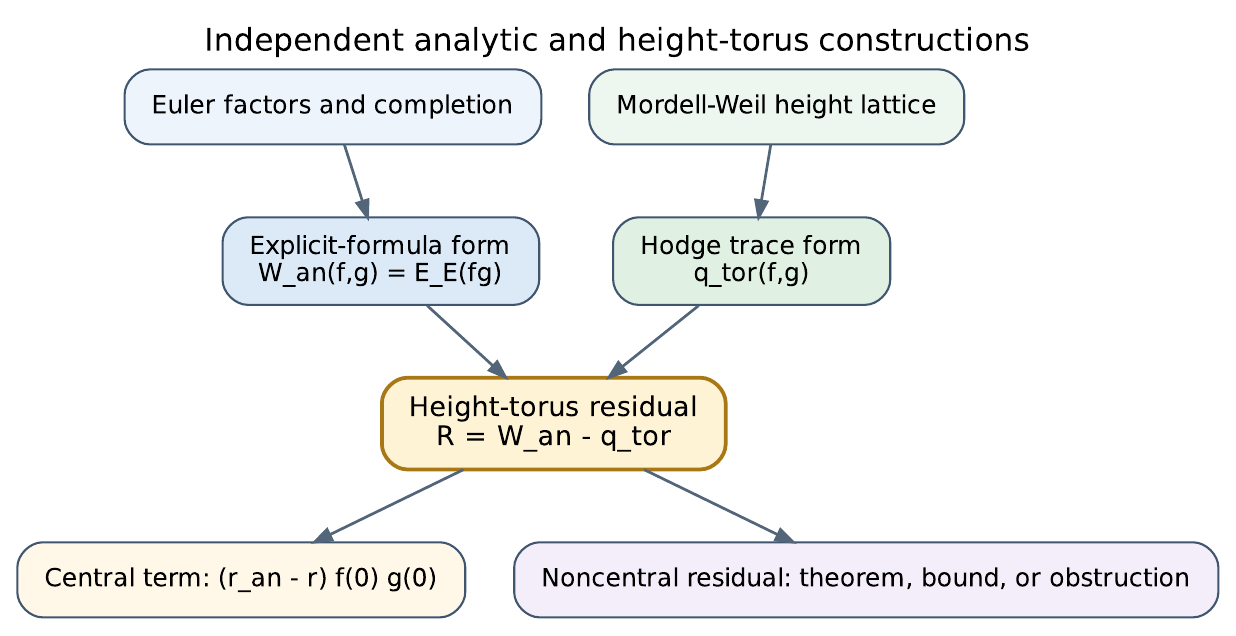}
  \caption{The analytic and metric-topological matrices are built by independent pipelines.  The residual may yield an exact finite match, a one-sided or norm comparison, a certified obstruction, or an inconclusive result requiring sharper error control.}
  \label{fig:dual-residual}
\end{figure}

\subsection{A counting-law obstruction to the raw torus spectrum}

The height torus records the algebraic rank in its zero modes, but its full positive spectrum has a counting law incompatible with the zero ordinates of a degree-two $L$-function.

\begin{theorem}[Raw flat-torus spectral no-go theorem]\label{thm:spectral-no-go}
Let $N_{\mathrm{tor}}(T)$ count the positive eigenvalues of $D_E^{(1)}$, with multiplicity, not exceeding $T$.  Let $N_{L(E)}(T)$ count the nontrivial zeros $\rho$ of $L(E,s)$ with $0<\operatorname{Im}\rho\le T$, with multiplicity.  Then
\begin{equation}\label{eq:counting-laws}
  N_{\mathrm{tor}}(T)=
  \begin{cases}
    0,&r=0,\\
    \Theta(T^r),&r\ge1,
  \end{cases}
  \qquad
  N_{L(E)}(T)=\Theta(T\log T).
\end{equation}
Therefore no affine rescaling of the raw spectrum of $D_E^{(1)}$, and no fixed finite change of multiplicity, can equal the multiset of all positive zero ordinates of $L(E,s)$.
\end{theorem}

\begin{proof}
For $r\ge1$, \cref{thm:flat-spectrum} gives positive frequencies $2\pi\norm{\xi}$ with $\xi\in M_E^*\setminus\{0\}$, each carrying the fixed one-form multiplicity $r$.  Lattice-point counting in an $r$-dimensional ellipsoid gives $N_{\mathrm{tor}}(T)=\Theta(T^r)$, equivalently the Weyl law for the compact flat torus \cite{Chavel1984}.  The standard Riemann--von Mangoldt theorem for a degree-two automorphic $L$-function gives $N_{L(E)}(T)=\Theta(T\log T)$ \cite{IwaniecKowalski2004}.  For $r=1$, the torus count lacks the factor $\log T$; for $r\ge2$, it has a different polynomial exponent; and for $r=0$ it is identically zero.  Affine rescaling and fixed multiplicity do not change these growth classes.
\end{proof}

\begin{figure}[H]
  \centering
  \includegraphics[width=0.82\textwidth]{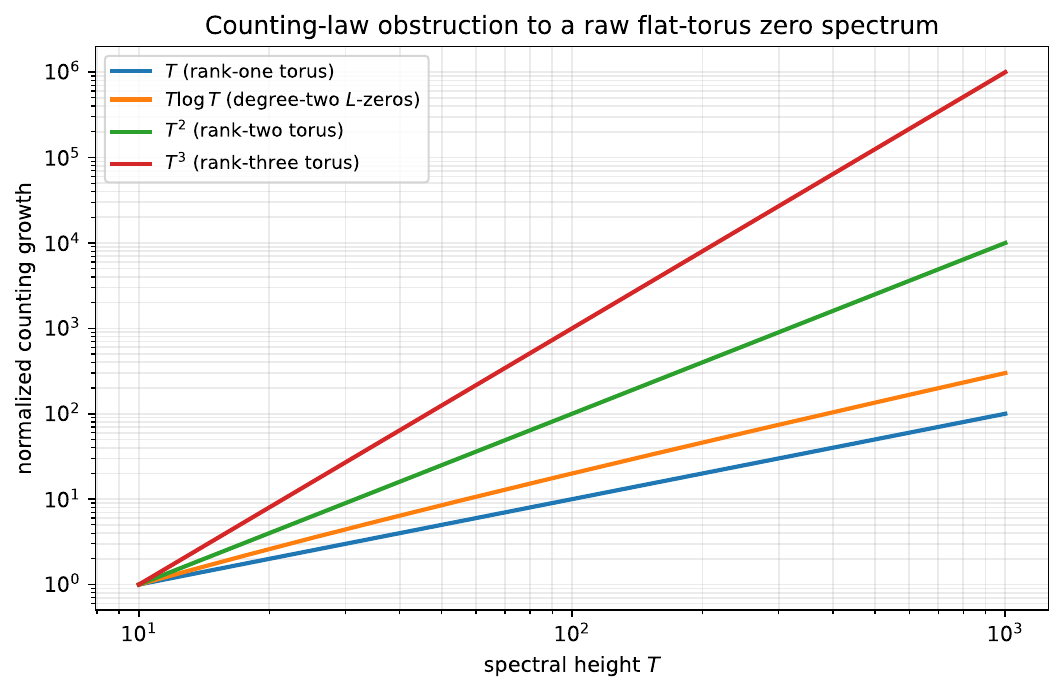}
  \caption{Structural comparison of counting laws, normalized at a common reference point.  The figure is not numerical zero data: it visualizes the proved incompatibility between $T\log T$ and the flat-torus laws $T^r$.}
  \label{fig:spectral-counting}
\end{figure}

\begin{corollary}[What the no-go theorem does not exclude]\label{cor:no-go-scope}
The theorem does not rule out the equality of central multiplicities
\[
  \ord_{s=1}L(E,s)=\dim\ker\Delta_E^{(1)},
\]
nor does it rule out a bridge containing an additional noncompact, cohomological, dynamical, or arithmetic operator.  It rules out only the claim that the raw compact height-torus spectrum, by a fixed linear rescaling, is already the complete zero spectrum.
\end{corollary}

\subsection{A phase cannot repair a Gram mismatch}

\begin{lemma}[Unitary phase invariance]\label{lem:unitary-phase}
Let $T:D(T)\subset H\to K$ be densely defined and let $U:K\to K$ be unitary.  Then, on $D(T)$,
\begin{equation}\label{eq:unitary-phase}
  \langle UTf,UTg\rangle_K=\langle Tf,Tg\rangle_K.
\end{equation}
Hence the quadratic forms generated by $UT$ and $T$ coincide.  If $T$ is closable (respectively closed), then $UT$ is closable (respectively closed), and their closed form closures and associated nonnegative self-adjoint operators coincide.
\end{lemma}

\begin{proof}
Unitarity gives $U^*U=I$, so
\[
  \langle UTf,UTg\rangle
  =\langle Tf,U^*UTg\rangle
  =\langle Tf,Tg\rangle.
\]
\end{proof}

\begin{corollary}[Scope of helical or Gamma phases]\label{cor:phase-scope}
A unimodular phase may impose a reality convention or rotate a representation, but it cannot change the torus Gram form, insert the real Gamma/digamma term in \eqref{eq:elliptic-explicit-formula}, or reduce a nonzero residual \eqref{eq:residual-form}.  The toroidal helices of \cref{sec:helices} therefore remain geometric representatives of homology classes, not phase-only solutions of the analytic bridge.
\end{corollary}

\section{Certified finite residual framework}\label{sec:residual-certification}

The residual construction is meaningful only if the two matrices are generated independently and the reported sign or norm survives all numerical errors.  This section specifies the finite certificate required before a residual observation is promoted to a theorem-level statement.

\subsection{Independent construction table}

\begin{table}[H]
\centering
\caption{The analytic and torus objects must be built independently before subtraction.}
\label{tab:independent-constructions}
\small
\begin{tabularx}{\textwidth}{@{}l>{\raggedright\arraybackslash}X>{\raggedright\arraybackslash}X>{\raggedright\arraybackslash}X@{}}
\toprule
Layer & Immutable input & Constructed object & Forbidden shared input \\
\midrule
Analytic & Euler factors, conductor, Gamma factor, fixed test convention & $G^{\mathrm{an}}_{E,\Delta}$ & Mordell--Weil basis, height Gram matrix, algebraic rank \\
Topological & Saturated Mordell--Weil basis, height normalization, dual lattice & $G^{\mathrm{tor}}_{E,\Delta,\tau}$ & Zeros, analytic rank, explicit-formula matrix \\
Comparison & Independently generated matrices and error bounds & $R_E=G^{\mathrm{an}}-G^{\mathrm{tor}}$ & A common factor used to manufacture both matrices \\
\bottomrule
\end{tabularx}
\end{table}

\subsection{Validated perturbation margins}

\begin{theorem}[Finite residual margin]\label{thm:residual-margin}
Let $R$ be the exact Hermitian residual matrix and let $\widetilde R$ be a computed Hermitian approximation satisfying
\begin{equation}\label{eq:residual-error}
  \norm{R-\widetilde R}_2\le\delta.
\end{equation}
Then every ordered eigenvalue satisfies
\begin{equation}\label{eq:residual-eigen-bound}
  \abs{\lambda_j(R)-\lambda_j(\widetilde R)}\le\delta.
\end{equation}
In particular:
\begin{enumerate}[label=\textup{(\roman*)}]
  \item if $\lambda_{\min}(\widetilde R)>\delta$, then $R\succ0$;
  \item if $\lambda_{\min}(\widetilde R)<-\delta$, then $R$ has a certified negative direction;
  \item if $\norm{\widetilde R}_2+\delta\le\varepsilon$, then $\norm{R}_2\le\varepsilon$.
\end{enumerate}
\end{theorem}

\begin{proof}
Equation \eqref{eq:residual-eigen-bound} is Weyl's Hermitian perturbation inequality.  The three conclusions follow by applying it to the smallest eigenvalue and by the triangle inequality for the operator norm; see \cite{Kato1976,Higham2002}.
\end{proof}

\begin{table}[H]
\centering
\caption{Permitted classifications for a finite residual calculation.}
\label{tab:residual-status}
\small
\begin{tabularx}{\textwidth}{@{}lXX@{}}
\toprule
Status & Validated condition & Permitted conclusion \\
\midrule
PASS--exact identity & Exact arithmetic proves $R=0$ entry by entry & Exact agreement on the stated finite test space only \\
PASS--controlled norm & $\norm{\widetilde R}_2+\delta\le\varepsilon$ for a predeclared $\varepsilon$ & Quantified approximation; not exact equality \\
PASS--one-sided & $\lambda_{\min}(\widetilde R)-\delta\ge0$ & $W^{\mathrm{an}}\ge q^{\mathrm{tor}}$ on that test space \\
FAIL--obstruction & $\lambda_{\min}(\widetilde R)+\delta<0$ & The tested domination or exact identity is impossible on that space \\
WARN & The validated interval intersects zero, or a tail/cross-term bound is missing & Numerical evidence only; no sign theorem \\
\bottomrule
\end{tabularx}
\end{table}

\subsection{Certificate and negative-control protocol}

\begin{protocol}[Auditable residual computation]\label{protocol:residual}
A residual run must record:
\begin{enumerate}
  \item the curve model, conductor, local-factor convention, Fourier convention, $\Delta$, $\tau$, and the Paley--Wiener seed family with its declared Fourier supports;
  \item immutable input hashes and software versions;
  \item a builder for $G^{\mathrm{an}}$ that includes conductor, Gamma, and all prime powers allowed by \cref{thm:exact-prime-support};
  \item a separate builder for $G^{\mathrm{tor}}$ from the saturated height Gram matrix and a certified dual-lattice enumeration;
  \item interval or outward-rounded bounds for the Archimedean integral, matrix assembly, and any omitted tails;
  \item the residual enclosure, the acceptance inequality, and a lightweight verifier that does not invoke the discovery code;
  \item at least the following negative controls: remove the Gamma term; change the Fourier normalization without translating the samples; replace a primitive generator by a nonprimitive multiple; perturb one local coefficient beyond the declared error; and apply a unitary phase while confirming \cref{lem:unitary-phase}.
\end{enumerate}
\end{protocol}

\begin{figure}[H]
  \centering
  \includegraphics[width=0.98\textwidth]{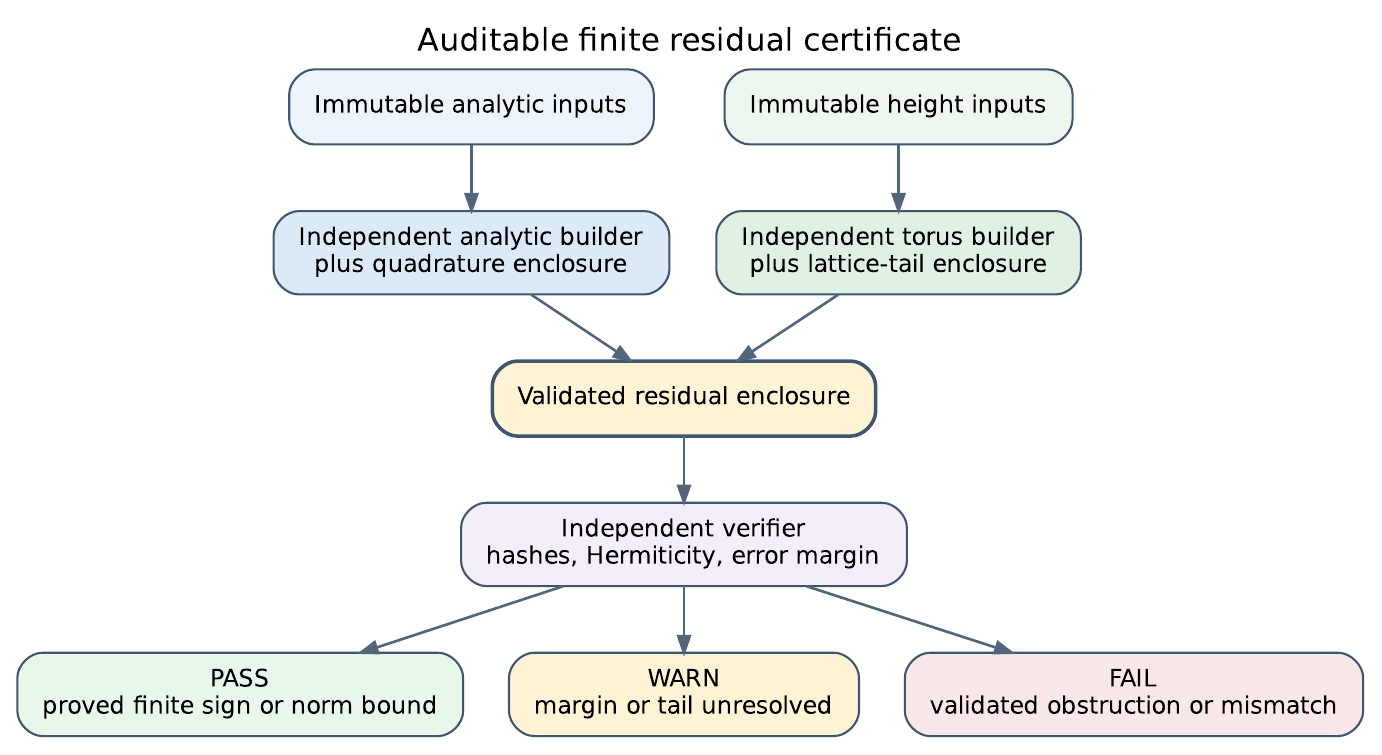}
  \caption{A residual certificate requires immutable inputs, independent matrix builders, a validated residual enclosure, and a verifier whose decision rule is fixed before the result is inspected.}
  \label{fig:residual-certificate}
\end{figure}

\subsection{Current implementation boundary}

The accompanying package implements and verifies the height-torus data, database-basis and bounded-saturation checks, Gram matrices, theta and heat calculations, and the structural figures.  It also supplies the exact definitions, matrix interfaces, support theorem, error inequalities, negative-control specification, and a deterministic residual-certificate verifier.  Synthetic PASS and FAIL records are included solely to test that the verifier accepts a positive margin and rejects a validated negative direction.  They are not elliptic-curve evidence.  The present manuscript does not report an unvalidated sign for $R_E$: a publishable residual table requires interval enclosure of the Archimedean integral and independent replay of both arithmetic builders.  This boundary is intentional.  A finite residual theorem or a certified obstruction is the next computational result to be earned, not a value to be inferred from a plot.

\section{BSD as a topological--spectral statement}\label{sec:bsd}

The preceding sections separate three claims that must not be conflated: the proved topology of the height torus, the independently defined explicit-formula functional, and a comparison residual.  The counting-law obstruction rules out a direct identification of the entire raw torus spectrum with all zero ordinates, but it does not decide the central multiplicity or the strong leading-coefficient formula.

\subsection{The rank conjecture}

Combining BSD with \cref{thm:rank-homology} yields an exact equivalence.

\begin{theorem}[Topological form of the BSD rank conjecture]\label{thm:topological-BSD}
For an elliptic curve $E/\Q$, the following statements are equivalent:
\begin{enumerate}[label=\textup{(\roman*)}]
  \item $\rank E(\Q)=\ord_{s=1}L(E,s)$;
  \item $b_1(\T_E)=\ord_{s=1}L(E,s)$;
  \item $\dim\ker\Delta_E^{(1)}=\ord_{s=1}L(E,s)$.
\end{enumerate}
\end{theorem}

\begin{proof}
The equality $b_1(\T_E)=\rank E(\Q)$ is \cref{thm:rank-homology}, and $\dim\ker\Delta_E^{(1)}=b_1(\T_E)$ follows from Hodge theory.  Substitution proves equivalence.
\end{proof}

\begin{warning}[Equivalence is not proof]\label{warn:not-proof}
The theorem moves the algebraic rank into topology but does not derive the order of vanishing from topology.  Any claimed proof must still establish the equality between the independently defined $L$-function and the independently defined harmonic geometry of $\T_E$.
\end{warning}

\subsection{The leading coefficient}

Using \cref{thm:reg-volume}, the full BSD formula \eqref{eq:BSD-full} becomes
\begin{equation}\label{eq:BSD-volume}
  L^*(E,1)
  =
  \frac{
    \Omega_E\,\Vol(\T_E)^2\,\#\Sha(E/\Q)\,\prod_pc_p
  }{
    \#E(\Q)_{\mathrm{tors}}^2
  }.
\end{equation}
Using the heat coefficient \eqref{eq:recover-volume}, it can also be written
\begin{equation}\label{eq:BSD-heat}
  L^*(E,1)
  =\frac{\Omega_E\,\#\Sha(E/\Q)\prod_pc_p}
  {\#E(\Q)_{\mathrm{tors}}^2}
  \left[
    \lim_{t\downarrow0}(4\pi t)^{r/2}\Tr(e^{-t\Delta_E})
  \right]^2.
\end{equation}
This is a metric-spectral reformulation of the regulator factor.  The remaining arithmetic factors are not absorbed into $\T_E$.

\subsection{Relation to regularized-determinant and trace-formula programs}

The demand for an independently defined bridge is not new.  Deninger's program seeks cohomological or dynamical objects whose infinitesimal generators express arithmetic local factors and zeta functions through regularized determinants \cite{Deninger1991,Deninger1992,Deninger2002}.  Connes' ad\`ele-class-space approach interprets the explicit formula as a noncommutative trace formula and gives a spectral realization of zeros in an absorption-spectrum framework \cite{Connes1999}.  In positive characteristic, Hesselholt constructed a topological-Hochschild-homology realization of the regularized-determinant picture for smooth proper schemes over finite fields \cite{Hesselholt2016}.

The present torus does something more limited and logically different.  Its Laplacian is built from the global Mordell--Weil lattice, so its harmonic multiplicity already equals the algebraic rank.  It neither reconstructs the Euler product nor supplies the missing arithmetic cohomology.  Accordingly, a future bridge may draw on Deninger- or Connes-type determinant and trace mechanisms, but it must also explain why their independently constructed arithmetic spectrum maps to the height-torus harmonic sector.  Merely declaring the zeros to be eigenvalues would not meet that requirement.

\subsection{The analytic bridge problem}

A valid bridge must be constructed without using the desired rank equality as an input.

\begin{problem}[Non-circular augmented bridge]\label{problem:bridge}
Construct an arithmetic, cohomological, dynamical, or noncommutative object $\mathscr A_E$, functorially from data of $E$ available independently of $\rank E(\Q)$, and prove both: (i) a determinant, trace, or cohomological comparison with the completed $L$-function; and (ii) a comparison with the height-torus harmonic sector.  Equivalently, construct an augmented torus-side form whose residual against $W^{\mathrm{an}}_{E,\Delta}$ admits a theorem-level identity, sign, norm bound, or finite-to-infinite limit, and whose central multiplicity equals $\dim\ker\Delta_E^{(1)}$ without inserting that dimension as input.
\end{problem}

At minimum, a proposed bridge must satisfy all of the following.

\begin{enumerate}[label=\textup{B\arabic*.}]
  \item \textbf{Independent definition.}  The construction must not insert $r$, a Mordell--Weil basis, or the zeros of $L(E,s)$ into the operator or boundary conditions whose spectrum is later said to recover them.
  \item \textbf{Functional analysis.}  If an operator is used, specify its Hilbert space, dense domain, closedness, self-adjointness or normality, boundary conditions, and spectral type.
  \item \textbf{Arithmetic comparison.}  Establish a theorem relating its determinant, trace, resolvent, or cohomology to the Euler product or completed $L$-function.
  \item \textbf{Central multiplicity.}  Prove, rather than assume, that the order at $s=1$ is the dimension of the relevant kernel or generalized eigenspace.
  \item \textbf{Leading coefficient and the $\Sha$ barrier.}  Account for the period, regulator, $\Sha$, Tamagawa, and torsion factors with correct normalizations.  In particular, the construction must prove the required finiteness and identify the order of $\Sha$; this is not a removable bookkeeping factor and may be as deep as the strong BSD conjecture itself.
  \item \textbf{Functoriality.}  Explain compatibility with isogenies, quadratic twists, base change, and local reduction data.
\end{enumerate}

\begin{warning}[The $\Sha$ factor is a separate open barrier]\label{warn:sha-barrier}
For ranks at least two, the general finiteness of $\Sha(E/\Q)$ is not known.  A determinant identity that simply leaves $\#\Sha$ as an unexplained multiplier has not completed the strong BSD comparison; a construction that claims to produce this factor must include an independent arithmetic mechanism for its finiteness and order.  The numerical quotients in \cref{sec:computations} are therefore reported as BSD consistency quotients, not as proofs of triviality of $\Sha$ in the rank-two and rank-three examples.
\end{warning}

The standard scalar Laplacian on $\T_E$ is useful because it rigorously packages the height lattice, but it is built from $M_E$ and hence already knows the algebraic rank.  It cannot by itself constitute an independent proof of BSD.  The value of the construction is that it identifies the exact target geometry a successful bridge would have to recover.  Moreover, \cref{thm:spectral-no-go} proves that its raw compact spectrum has the wrong high-frequency counting law to equal the full set of zero ordinates.  A successful bridge must therefore add arithmetic structure rather than merely rename the flat-torus Laplacian.

\section{Computational realization and verified examples}\label{sec:computations}

The computations in this section were generated with PARI/GP 2.17.2 through \texttt{cypari2}, using Cremona/LMFDB curve data where available \cite{Cremona1997,LMFDB2016,PARI2026,CremonaEcdata2026}.  The implementation follows standard computational frameworks for canonical heights and motivic $L$-values \cite{CremonaPrickettSiksek2006,Dokchitser2004}.  The complete script, numerical JSON data, CSV summary, and vector figures accompany the manuscript.

\subsection{Protocol and normalization}

\begin{protocol}[Reproducible height-torus computation]\label{protocol:computation}
For each curve:
\begin{enumerate}
  \item initialize the general Weierstrass model and compute the conductor, root number, torsion structure, and rank bounds;
  \item retrieve the database $\Z$-basis with \texttt{ellgenerators} or \texttt{ellidentify} whenever the curve is covered by the installed Cremona \texttt{elldata} package;
  \item treat points returned by \texttt{ellrank} only as witnesses or candidate generators, since they need not be saturated; compare any supplied candidates with the database basis and use \texttt{ellsaturation} as an additional bounded-prime check;
  \item compute the height Gram matrix with \texttt{ellheightmatrix}; for a full basis its determinant is the regulator in the Cremona/PARI normalization;
  \item compute the volume, systole, and successive minima by exact lattice formulas followed by finite enumeration justified by the least Gram eigenvalue;
  \item compute the analytic rank and the first nonzero derivative with PARI's $L$-function routines, and report the resulting BSD consistency quotient without promoting it to an unconditional computation of $\Sha$.
\end{enumerate}
\end{protocol}

All displayed Gram matrices, squared lengths, regulators, volumes, and BSD consistency quotients use the direct Cremona/PARI height $\qheight$.  Equation \eqref{eq:cp-silverman-conversion} gives the optional conversion to the Silverman convention; no conversion is performed inside the numerical BSD check.

\subsection{Curves and rational points}

\begin{figure}[H]
  \centering
  \includegraphics[width=0.95\textwidth]{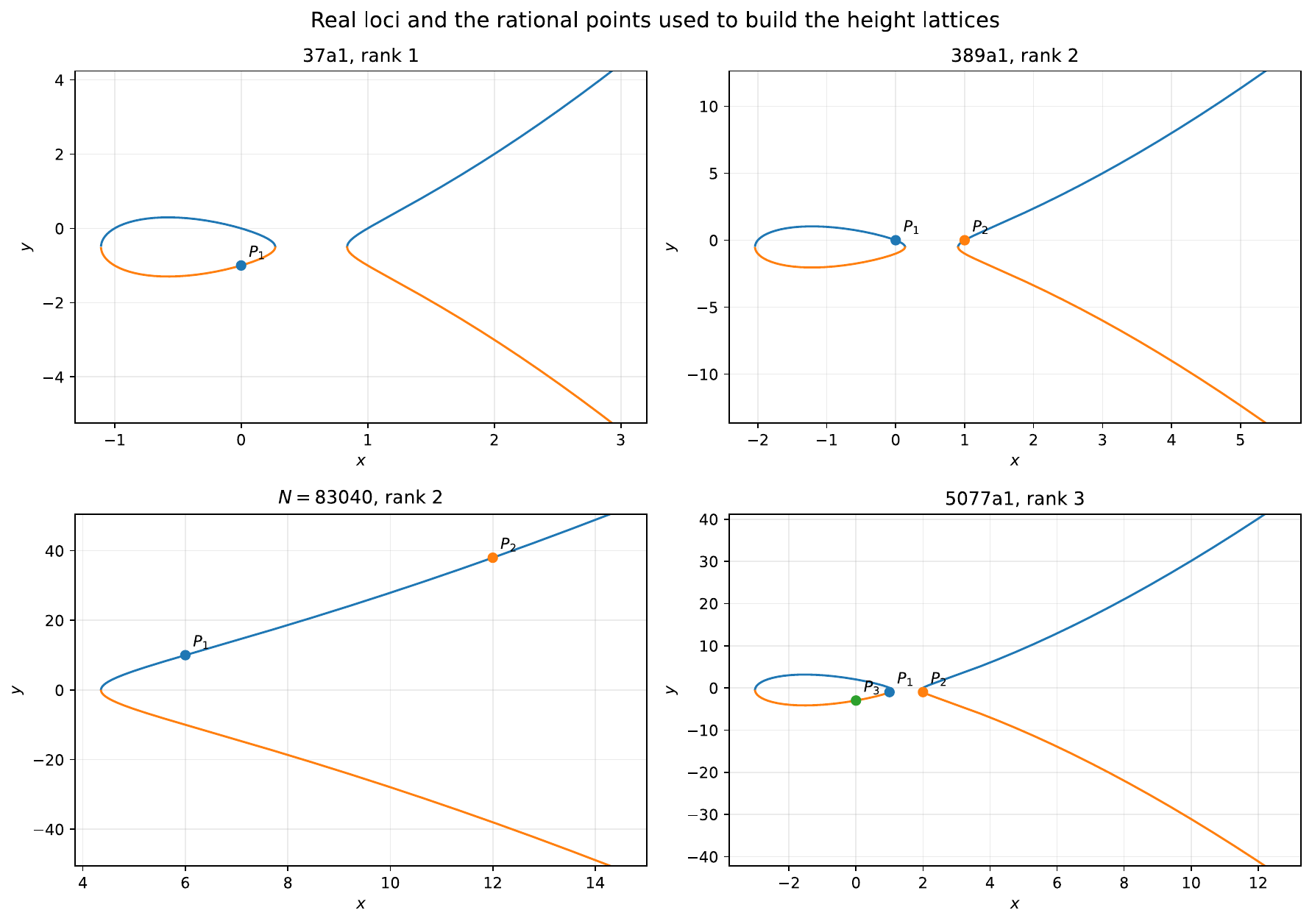}
  \caption{Real loci and rational points used to construct the example height lattices.  The geometric positions in the plane do not encode independence; independence is measured by the N\'eron--Tate Gram matrix.}
  \label{fig:curves-generators}
\end{figure}

\begin{table}[H]
\centering
\caption{Computed arithmetic and height-torus data.  Values are rounded.  Every positive-rank row uses a full $\Z$-basis from the Cremona/PARI database; supplied candidate points were compared with that basis and have index one.  Bounded calls to \texttt{ellsaturation} are retained as independent consistency checks, not as substitutes for the database basis certificate.}
\label{tab:curve-summary}
\small
\begin{tabularx}{\textwidth}{@{}l>{\raggedright\arraybackslash}Xrrrrr@{}}
\toprule
Curve & Weierstrass equation & $N$ & $r$ & $r_{\rm an}$ & $w$ & $\#E(\Q)_{\rm tors}$ \\
\midrule
11a1 & $y^2+y=x^3-x^2-10x-20$ & 11 & 0 & 0 & $+1$ & 5 \\
37a1 & $y^2+y=x^3-x$ & 37 & 1 & 1 & $-1$ & 1 \\
389a1 & $y^2+y=x^3+x^2-2x$ & 389 & 2 & 2 & $+1$ & 1 \\
$N=83040$ & $y^2=x^3-x^2-10x-20$ & 83040 & 2 & 2 & $+1$ & 1 \\
5077a1 & $y^2+y=x^3-7x+6$ & 5077 & 3 & 3 & $-1$ & 1 \\
\bottomrule
\end{tabularx}

\medskip
\begin{tabularx}{0.92\textwidth}{@{}lrrr>{\raggedright\arraybackslash}X@{}}
\toprule
Curve & $\Reg$ & $\Vol$ & $\sys^2$ & basis status \\
\midrule
11a1 & 1.000000 & 1.000000 & -- & rank-zero database case \\
37a1 & 0.051111 & 0.226078 & 0.051111 & database basis; index $1$ \\
389a1 & 0.152460 & 0.390461 & 0.327001 & database basis; index $1$ \\
$N=83040$ & 3.887410 & 1.971652 & 1.708471 & database basis; index $1$ \\
5077a1 & 0.417144 & 0.645867 & 0.668205 & database basis; index $1$ \\
\bottomrule
\end{tabularx}
\end{table}

The bases or subgroup generators used are
\begin{align*}
  37\mathrm{a}1 &: P_1=(0,-1),\\
  389\mathrm{a}1 &: P_1=(0,0),\quad P_2=(1,0),\\
  N=83040 &: P_1=(6,10),\quad P_2=(12,38),\\
  5077\mathrm{a}1 &: P_1=(1,-1),\quad P_2=(2,-1),\quad P_3=(0,-3).
\end{align*}
Every listed point was checked directly against its Weierstrass equation.  For the rank-two curve from the original project, this corrects earlier plotted points that did not satisfy the equation.  In the clean PARI/GP 2.17.2 verification, the displayed positive-rank point sets agree with full database Mordell--Weil bases up to a unimodular change of basis; equivalently, the regulator ratio gives index $1$.  Calls to \texttt{ellsaturation} leave these bases unchanged at the stated check bounds.  No unsaturated witness set returned by \texttt{ellrank} is used to form a regulator.

\subsection{Basis, saturation, and BSD consistency checks}

For a candidate basis $P$ and a database basis $G$, the index is checked from the covolume identity
\[
  [\langle G\rangle:\langle P\rangle]
  =\sqrt{\frac{\det\operatorname{Gram}(P)}{\det\operatorname{Gram}(G)}}.
\]
The clean run also evaluates
\begin{equation}\label{eq:analytic-sha-quotient}
  S_{\mathrm{an}}(E)
  :=\frac{L^{(r)}(E,1)/r!}{\operatorname{ellbsd}(E)\Reg(E/\Q)}.
\end{equation}
PARI's definition predicts $S_{\mathrm{an}}(E)=\#\Sha(E/\Q)$ if the full BSD formula holds in the same normalization.  The following table records the reproducibility checks; $S_{\mathrm{an}}\approx1$ is a numerical BSD-II consistency result, not an unconditional proof of $\Sha=0$ in ranks two and three.

\begin{table}[H]
\centering
\caption{Independent basis and BSD consistency checks in PARI/GP 2.17.2.  The saturation bound means that no prime below that bound divides the remaining subgroup index; global index one is supplied by agreement with the database $\Z$-basis.}
\label{tab:basis-bsd-verification}
\small
\begin{tabular}{@{}lccc@{}}
\toprule
Curve & saturation check $B$ & index to database basis & $S_{\mathrm{an}}(E)$ \\
\midrule
11a1 & -- & -- & 1.000000 \\
37a1 & 100 & 1 & 1.000000 \\
389a1 & 100 & 1 & 1.000000 \\
$N=83040$ & 5000 & 1 & 1.000000 \\
5077a1 & 100 & 1 & 1.000000 \\
\bottomrule
\end{tabular}
\end{table}

\subsection{Height Gram matrices}

For $37\mathrm{a}1$,
\[
  G_{37\mathrm{a}1}
  =\begin{pmatrix}0.05111140824\end{pmatrix}.
\]
For $389\mathrm{a}1$,
\begin{equation}\label{eq:G389}
  G_{389\mathrm{a}1}
  =\begin{pmatrix}
  0.3270007737 & 0.0585226748\\
  0.0585226748 & 0.4767116593
  \end{pmatrix}.
\end{equation}
For the full database basis on the curve of conductor $83040$, independently unchanged by the $5000$-saturation check,
\begin{equation}\label{eq:G83040}
  G_{83040}
  =\begin{pmatrix}
  1.7084707350 & 0.0491898454\\
  0.0491898454 & 2.2767905493
  \end{pmatrix}.
\end{equation}
For $5077\mathrm{a}1$,
\begin{equation}\label{eq:G5077}
  G_{5077\mathrm{a}1}
  =\begin{pmatrix}
  0.6682051657 & 0.0333380078 & -0.2365919007\\
  0.0333380078 & 0.7670433553 & -0.2764342921\\
  -0.2365919007 & -0.2764342921 & 0.9909063332
  \end{pmatrix}.
\end{equation}
The determinants agree with the regulators reported in \cref{tab:curve-summary}.  The positive eigenvalues of each matrix verify positive definiteness numerically; database and saturation information provide the separate arithmetic basis status.

\begin{figure}[H]
  \centering
  \includegraphics[width=0.86\textwidth]{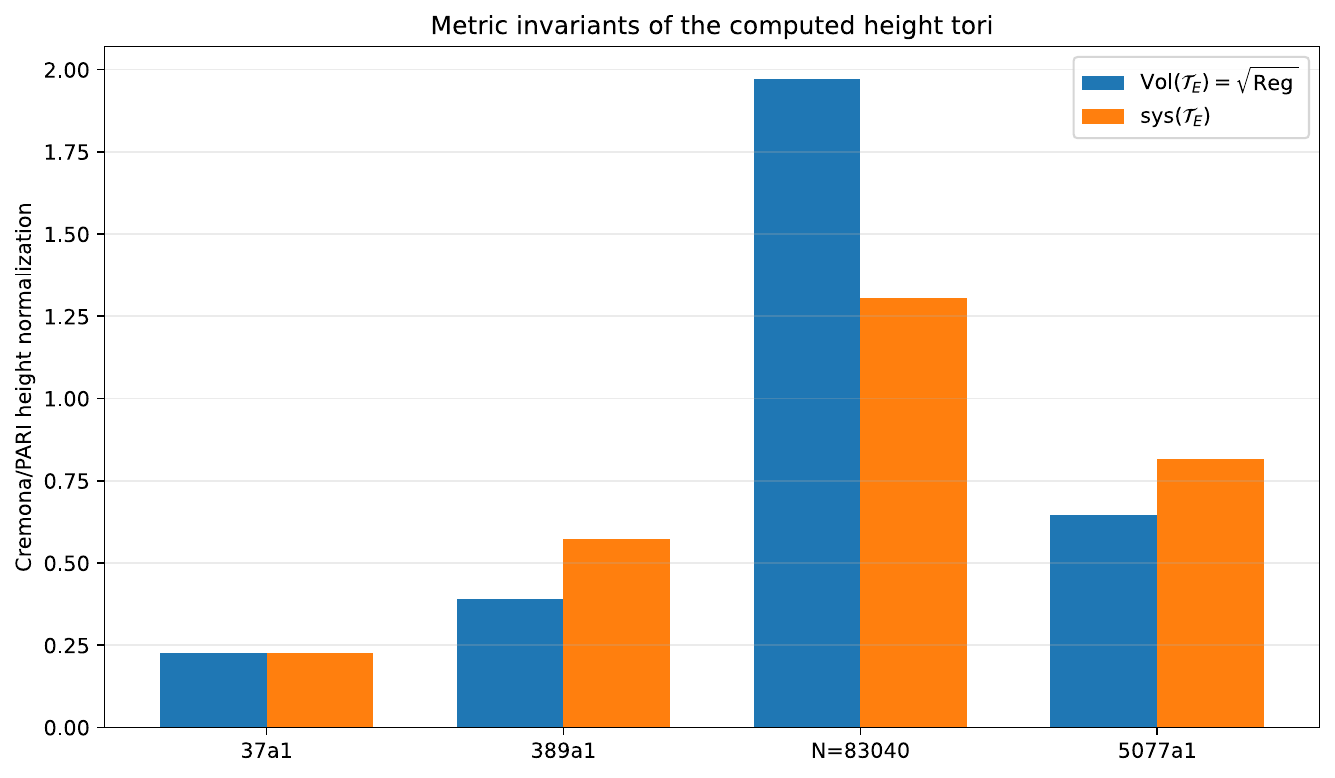}
  \caption{Volume and systole in the Cremona/PARI normalization.  Rank is a topological invariant; volume and systole discriminate the metric geometry among tori of the same or different rank.}
  \label{fig:volume-systole}
\end{figure}

\subsection{Central \texorpdfstring{$L$}{L}-function behavior}

\begin{figure}[H]
  \centering
  \includegraphics[width=0.82\textwidth]{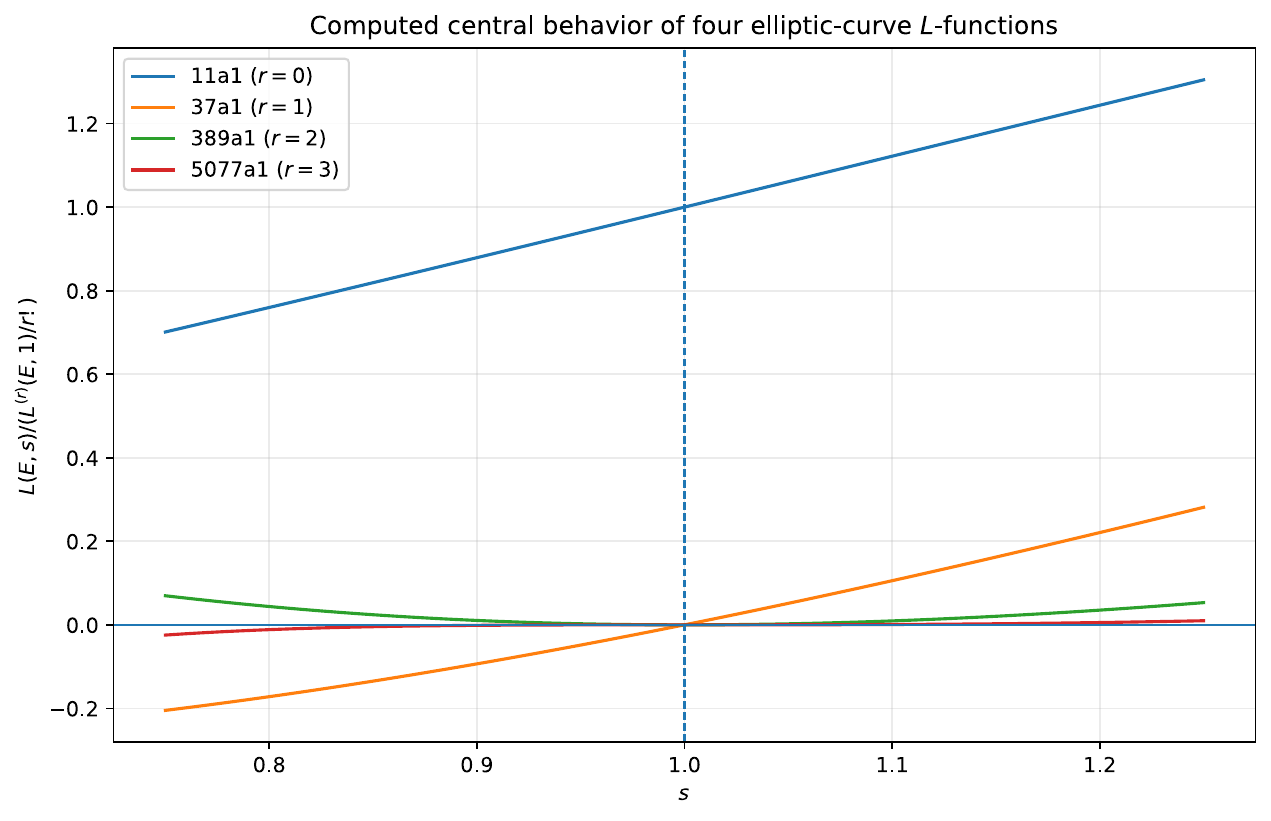}
  \caption{Numerically computed values of four elliptic-curve $L$-functions near $s=1$, divided by $L^{(r)}(E,1)/r!$; the curves have analytic ranks $0,1,2,3$.  Unlike a schematic plot of $(s-1)^r$, these curves were evaluated by PARI's $L$-function routines; the figure illustrates rather than proves the orders of vanishing.}
  \label{fig:L-central}
\end{figure}

For all five rows, the computed leading coefficient and the PARI BSD factor satisfy
\[
  S_{\mathrm{an}}(E)
  =\frac{L^*(E,1)}{\operatorname{ellbsd}(E)\Reg(E)}
  \approx1
\]
at the displayed precision, using full database bases and the direct Cremona/PARI regulator.  Thus the strong-BSD formula predicts $\#\Sha(E/\Q)=1$ for these examples.  For the rank-two and rank-three curves this sentence is deliberately conditional: the numerical quotient neither proves finiteness nor independently proves triviality of $\Sha$; see \cref{warn:sha-barrier}.

\subsection{Data integrity and reproducibility}

The figures in this paper fall into two classes.

\begin{itemize}
  \item \textbf{Structural figures} such as \cref{fig:pipeline,fig:rank-atlas,fig:helix,fig:4d-windows,fig:exact-support,fig:dual-residual,fig:spectral-counting,fig:residual-certificate} illustrate proved constructions, asymptotic classes, or verification architecture and are not presented as empirical evidence.
  \item \textbf{Computed figures} such as \cref{fig:curves-generators,fig:rank2-lattices,fig:theta,fig:heat,fig:volume-systole,fig:L-central} are generated from the accompanying arithmetic data and script.
\end{itemize}

This separation prevents a drawn equality from being mistaken for a numerical test.  In particular, no plot in the final manuscript is generated by setting the vertical data equal to the expected rank in advance.

\section{Prime-sum diagnostics and their proper role}

The Euler coefficients $a_p$ are the local data from which $L(E,s)$ is built, so prime sums remain natural for exploratory computation.  The logical requirements are:

\begin{enumerate}
  \item the statistic must respect the Hasse bound and its actual normalization;
  \item the rank must not be supplied to the feature construction;
  \item training and test sets must be separated by conductor or another meaningful arithmetic regime;
  \item comparisons must include established Mestre--Nagao baselines;
  \item finite-range correlations must not be promoted to asymptotic theorems without proof.
\end{enumerate}

A representative Mestre--Nagao-type sum is
\begin{equation}\label{eq:Mestre-Nagao}
  S(E,B)=\frac1{\log B}
  \sum_{\substack{p<B\\p\nmid N_E}}
  \frac{a_p\log p}{p},
\end{equation}
with several variants in the literature.  Such sums arise from explicit-formula heuristics and have a theoretical relation to average ranks in elliptic surfaces \cite{Nagao1997,RosenSilverman1998}.  Modern computational studies show that performance depends on the cutoff, conductor range, combination of multiple sums, and learned weights \cite{KazalickiVlah2022,BujanovicKazalickiNovak2024,BujanovicKazalickiVlah2025,BieriEtAl2026}.  Therefore any future computational extension of the present framework should treat prime sums as independent diagnostics for the analytic side, not as substitutes for the proved height-torus geometry.

A potentially useful research design is to compare two independently generated objects:
\[
  \text{local analytic features from }(a_p)
  \quad\text{versus}\quad
  \text{global metric features from }G_E.
\]
The comparison is scientifically meaningful only when the height features are withheld from the analytic predictor and the analytic rank is not inserted into the topological construction.

\section{Generalization to number fields and abelian varieties}

\subsection{Number fields}

Let $K$ be a number field.  The Mordell--Weil theorem gives
\[
  M_{E,K}:=E(K)/E(K)_{\mathrm{tors}}\cong\Z^{r_K}.
\]
The global N\'eron--Tate height defines a positive-definite form on $M_{E,K}\otimes\R$, and one may set
\begin{equation}\label{eq:number-field-torus}
  \T_{E,K}:=(M_{E,K}\otimes\R)/M_{E,K}.
\end{equation}
Then
\[
  H_1(\T_{E,K},\Z)\cong M_{E,K},
  \qquad
  b_1(\T_{E,K})=\rank E(K),
\]
and the volume is the square root of the regulator in the chosen global normalization.  The full BSD formula over $K$ contains a product of archimedean periods and local factors; the torus again accounts for the free Mordell--Weil rank and regulator, not the entire formula.

\subsection{Polarized abelian varieties}

For an abelian variety $A/K$, a symmetric ample line bundle or polarization $\lambda$ defines a N\'eron--Tate quadratic form on $A(K)/A(K)_{\mathrm{tors}}$ \cite{Lang1983,HindrySilverman2000,BombieriGubler2006,Pazuki2015}.  Define
\begin{equation}\label{eq:abelian-height-torus}
  \T_{A,\lambda}
  :=\Bigl(\bigl(A(K)/A(K)_{\mathrm{tors}}\bigr)\otimes_{\Z}\R\Bigr)
    \Big/\bigl(A(K)/A(K)_{\mathrm{tors}}\bigr).
\end{equation}
The same homology, length, volume, theta, and heat statements apply to the paired Mordell--Weil lattice.  In higher dimension, however, the polarization is part of the data, and the regulator in the general BSD formula is naturally a pairing between $A(K)$ and the dual abelian variety.  Thus the elliptic-curve case is the cleanest self-dual model, not a justification for ignoring polarization choices in general.

\section{Relation to existing work, novelty, and limitations}

\subsection{What is established mathematics}

The following ingredients are classical or standard:

\begin{itemize}
  \item finite generation of $E(\Q)$ \cite{Mordell1922,Weil1929};
  \item canonical heights and the regulator \cite{Neron1965,Silverman2009,Silverman1994};
  \item the interpretation of the free group as a Euclidean lattice \cite{Shioda1990,Shioda1991,SchuettShioda2019};
  \item lattice covolume, successive minima, and theta transformation \cite{Minkowski1910,Cassels1959,ConwaySloane1999};
  \item Fourier spectrum and heat kernel of a flat torus \cite{BergerGauduchonMazet1971,Chavel1984,McKeanSinger1967};
  \item regularized-determinant, dynamical, and noncommutative trace-formula programs for arithmetic zeta functions \cite{Deninger1991,Deninger1992,Deninger2002,Connes1999,Hesselholt2016}.
\end{itemize}

These results are included with references because they strengthen the framework; their prior existence does not invalidate their use.  It does prevent them from being presented as newly proved arithmetic theorems.

\subsection{Contribution of the present paper}

The contribution is the following integrated and corrected BSD-specific formulation:

\begin{enumerate}
  \item an explicit obstruction analysis showing why the topology of $E(\R)$ and a countable multiple orbit cannot carry arbitrary arithmetic rank;
  \item a single canonical object $\T_E$ in which rank, height, regulator, systole, and heat spectrum become different levels of one invariant dictionary;
  \item the marked-length reconstruction formula \eqref{eq:cross-reconstruction} as a finite geometric recovery procedure for the height matrix;
  \item the finite-index covering interpretation of unsaturated Mordell--Weil subgroups;
  \item a rigorous reconstruction of the helical idea through winding classes on $\T_E$, together with the four-dimensional character-window theorem;
  \item a centered explicit-formula functional and a Paley--Wiener seed construction that produce an analytic matrix independently of the height lattice while aligning the analytic and torus zero modes;
  \item the exact support theorem \cref{thm:exact-prime-support}, which converts compact Fourier support into a genuinely finite prime-power side;
  \item the height-torus residual \eqref{eq:residual-form}, the central-mode decomposition \eqref{eq:central-mode-decomposition}, congruence covariance, and the finite certificate framework of \cref{sec:residual-certification};
  \item the spectral-counting no-go theorem \cref{thm:spectral-no-go}, which proves that the raw compact torus cannot already be the complete zero operator;
  \item exact topological and heat-coefficient forms of the two BSD statements, accompanied by explicit non-circularity and $\Sha$ barriers;
  \item a reproducible computational realization that distinguishes certified bases, bounded saturation, schematic graphics, actual $L$-function calculations, and future residual claims requiring independent validation.
\end{enumerate}

The standard ingredients are cited at their original sources.  The intended contribution is the corrected BSD-specific integration, the explicit residual problem, the obstruction theorem for the raw spectrum, and a falsifiable route to the next comparison result.

\subsection{Limitations}

The framework has clear boundaries.

\begin{itemize}
  \item Constructing $\T_E$ requires knowledge of the free Mordell--Weil group or at least a full-rank subgroup.  It is therefore not, by itself, an algorithm for discovering the rank from local data.
  \item The torus omits torsion and the local/global obstruction group $\Sha$.  In ranks at least two, general finiteness of $\Sha$ is itself open, so recovering it requires substantially more than adding a formal decoration or an unexplained determinant factor.
  \item The scalar Laplace spectrum can have isospectral ambiguity.  The arithmetic marking contains more information than the unmarked spectrum.
  \item A toroidal helix is a visualization of an arithmetic class, not evidence that the Hasse--Weil $L$-function is a contour integral over that helix.
  \item The computational examples verify the internal dictionary; they do not test BSD independently when both algebraic and analytic ranks are imported from established software or databases.
  \item The scale $\tau$ in the torus trace form is declared rather than inferred; a bridge theorem must explain its arithmetic normalization or replace it with an intrinsic construction.
  \item The central-mode decomposition isolates $r_{\mathrm{an}}-r$, but the noncentral residual remains uncontrolled; the decomposition is not a proof that the rank discrepancy vanishes.
  \item The finite residual framework does not by itself provide an infinite-dimensional comparison.  Trace, determinant, or form limits require the topology and tail estimates appropriate to the claimed conclusion.
  \item The analytic bridge in \cref{problem:bridge} remains open.
\end{itemize}

\section{Conclusion}

The topological perspective becomes rigorous once the topology is attached to the Mordell--Weil lattice rather than to a bounded drawing of the real elliptic curve.  The resulting height torus
\[
  \T_E=(E(\Q)/E(\Q)_{\mathrm{tors}}\otimes\R)
  /(E(\Q)/E(\Q)_{\mathrm{tors}})
\]
has an exact arithmetic interpretation:
\begin{equation}\label{eq:final-dictionary}
  \boxed{
  \begin{aligned}
    H_1(\T_E,\Z)&\cong E(\Q)/E(\Q)_{\mathrm{tors}},\\
    b_1(\T_E)&=\rank E(\Q),\\
    \operatorname{Length}(\gamma_P)^2&=\qheight(P),\\
    \Vol(\T_E)^2&=\Reg(E/\Q),\\
    \Spec(\Delta_E)&=\{4\pi^2\norm{\xi}^2:\xi\in M_E^*\}.
  \end{aligned}}
\end{equation}
The original intuition that independent rational generators should correspond to independent loops is therefore correct after replacing the discrete orbit by the genuine geodesic $tP\bmod M_E$ and replacing the fixed topology of $E(\R)$ by the rank-dimensional quotient of its arithmetic lattice.

This yields a precise topological version of BSD:
\[
  \ord_{s=1}L(E,s)=b_1(\T_E),
\]
and a metric version of the leading coefficient in which the regulator is the squared volume of $\T_E$.  The paper does not claim that this equivalence proves BSD.

The analytic development fixes a second, independent object: the explicit-formula form $W^{\mathrm{an}}_{E,\Delta}$.  Its comparison with the positive torus trace form is encoded by the residual $\mathcal R_{E,\Delta,\tau}$.  The Paley--Wiener seed map aligns the zero-frequency terms, giving the exact decomposition
\[
  \mathcal R_{E,\Delta,\tau}(f,g)
  =\bigl(r_{\mathrm{an}}(E)-r\bigr)f(0)g(0)
   +\mathcal R^{\circ}_{E,\Delta,\tau}(f,g).
\]
The exact support theorem makes the local side finitely auditable, while the spectral-counting theorem proves that the raw compact torus cannot by itself be the complete zero operator.  Thus the paper identifies a precise comparison object, exposes the central rank discrepancy without assuming its vanishing, proves a structural obstruction for the simplest spectral identification, and specifies the certificate required for a valid finite result.

A successful continuation must therefore add an independently motivated arithmetic layer and prove what happens to the residual, the central multiplicity, the limiting topology, and the strong BSD factors.  Operator constructions, trace identities, arithmetic cohomology, and noncommutative or dynamical models can now be tested against fixed analytic and metric-topological targets rather than against a visual analogy.

\appendix
\section{Computational formulas and commands}

For a curve initialized in PARI/GP as \texttt{E=ellinit([a1,a2,a3,a4,a6])}, the core operations used were conceptually:
\begin{verbatim}
ellglobalred(E)         /* conductor and global reduction data */
ellrootno(E)            /* global root number */
elltors(E)              /* rational torsion */
ellrank(E, effort)      /* rank bounds and witness/candidate points */
ellgenerators(E)        /* database Z-basis, when elldata covers E */
ellidentify(E)          /* label, minimal model, database basis, change map */
ellsaturation(E, P, B)  /* remove index primes below B only */
ellheightmatrix(E, P)   /* Cremona/PARI Neron-Tate Gram matrix */
ellanalyticrank(E)      /* analytic rank and leading derivative data */
ellL1(E, r)             /* r-th derivative at the central point */
ellbsd(E)               /* c with L^(r)(1)/r! = c * Reg * #Sha predicted */
\end{verbatim}
The accompanying Python program performs the same calculations through \texttt{cypari2}, exports JSON and CSV data, and generates all vector figures.  A separate plain-GP script, \path{verify_basis_and_bsd.gp}, retrieves database bases, compares regulator covolumes, runs bounded saturation checks, and evaluates \eqref{eq:analytic-sha-quotient}.  The theta and heat plots use finite symmetric boxes in $\Z^r$; they are visual approximations to absolutely convergent lattice sums, not symbolic evaluations of the infinite series.

\paragraph{Important basis warning.}
The point vector returned by \texttt{ellrank} is evidence for the lower rank bound, but it is not promised to be a saturated basis.  The PARI manual's own rank-three example shows a witness set whose determinant drops by a factor of $9$ after saturation.  Therefore, the final verification script uses \texttt{ellgenerators}/\texttt{ellidentify} as the global basis certificate when database coverage is available and uses \texttt{ellsaturation} only as a bounded-prime consistency test.

\section{A basis-status checklist}

Before a numerical Gram determinant is called the regulator of $E/\Q$, the following must be recorded:

\begin{enumerate}
  \item proof or database source for the algebraic rank;
  \item exact coordinates of the proposed generators and direct verification that they lie on the curve;
  \item nonsingularity and minimal-model information;
  \item linear independence under the canonical height pairing;
  \item saturation status, including the bound and whether a global basis theorem is available;
  \item normalization of the height matrix;
  \item numerical precision and stability of the determinant.
\end{enumerate}

A positive Gram determinant proves real linear independence of the supplied classes, not that their span has index one.  The latter is a separate arithmetic assertion measured geometrically by \cref{thm:index-cover}.

\section{Transform and residual convention sheet}

The residual construction uses the Fourier convention \eqref{eq:residual-fourier}; the completed function \eqref{eq:completed-E-new}; the generalized local coefficients \eqref{eq:Lambda-E}; the real even Paley--Wiener seed space \eqref{eq:test-seed-space} with $h_{f,g}=fg$ and \eqref{eq:seed-convolution}; and an explicitly declared torus scale $\tau$.  A residual computation is invalid if any of these conventions changes between the analytic and topological builders without translating every sample, multiplier, and normalization.

For a finite seed family $\mathcal F$, the minimum machine-readable record is
\begin{center}
\small
\begin{tabularx}{0.94\textwidth}{@{}lX@{}}
\toprule
Field & Required content \\
\midrule
curve & Weierstrass coefficients, label, conductor, bad-prime data \\
analytic convention & completion, Fourier convention, $\Delta$, basis functions, Gamma quadrature rule \\
height convention & Cremona/PARI or converted convention, saturated basis source, Gram matrix \\
torus convention & $\tau$, dual-lattice enumeration bound, omitted-tail enclosure \\
certificate & residual matrix or interval matrix, norm bound $\delta$, acceptance inequality \\
provenance & software versions, precision, immutable hashes, independent verifier output \\
\bottomrule
\end{tabularx}
\end{center}

\section{Residual negative controls}

Every confirmatory run must include at least one algebraic and one numerical failure case.  The expected results are:
\begin{enumerate}
  \item omitting the Gamma term changes $G^{\mathrm{an}}$ beyond the validated assembly error;
  \item changing the Fourier convention while retaining the old sample locations fails the basis-covariance or replay check;
  \item replacing a primitive generator by $mP$ changes the candidate regulator by an index-square factor;
  \item a corrupted local coefficient is detected by the immutable-input hash or by independent reconstruction;
  \item left multiplication by a unitary phase leaves the constructed Gram form unchanged, as required by \cref{lem:unitary-phase};
  \item declaring an infinite conclusion without a tail or complement theorem is classified as WARN rather than PASS.
\end{enumerate}

\section*{Statements and declarations}

\noindent\textbf{Competing interests.} The author declares no competing interests.

\medskip
\noindent\textbf{Data and code availability.} The complete LaTeX source, bibliography, computational scripts, plain-PARI verification script, verification summary, generated arithmetic data, structural residual-figure generator, and vector figures are supplied with the manuscript.  The package distinguishes the verified height-torus computations from the residual certificate protocol, for which no unvalidated sign is reported.

\medskip
\noindent\textbf{Research status.} The manuscript gives a rigorous reformulation and a computational framework.  It does not claim a proof of the Birch and Swinnerton--Dyer conjecture.

\bibliographystyle{abbrvnat}
\bibliography{references}

\end{document}